\documentclass{amsart}
\usepackage[english]{babel}
\usepackage[utf8]{inputenc}
\usepackage{graphicx}
\usepackage[colorinlistoftodos]{todonotes}
\usepackage{amsmath}
\usepackage{amsfonts,amssymb,amsthm}
\usepackage{txfonts}
\usepackage{verbatim}
\usepackage[active]{srcltx} 
\usepackage{enumerate} 
\usepackage{relsize}

\usepackage{url}
\newcommand{\To}{\longrightarrow}
\newcommand{\lra}{\longrightarrow} \newcommand{\ra}{\rightarrow}
\newcommand{\h}{\mathcal{H}}

\newcommand{\D}{\mathcal{D}}
\newcommand{\s}{\boldsymbol{\mathcal{S}}}
\newcommand{\A}{\mathcal{A}}
\newcommand{\R}{\mathbb{R}}

\newcommand{\K}{\mathcal{K}}

\newcommand{\M}{\mathcal{M}}

\newcommand{\W}{\mathcal{W}} 
\newcommand{\bW}{\mathbf{W}}
\newcommand{\G}{\mathcal{G}}

\newcommand{\spec}{\mathrm{spec}}
\newcommand{\N}{\mathbb{N}}

\newcommand{\je}{J_{\Elk}}

\renewcommand{\ker}{\operatorname{ker}}
\newcommand{\Rang}{\operatorname{Rang}}
\newcommand{\Real}{\mathbb{R}}

\newcommand{\Complex}{\mathbb{C}}

\newcommand{\F}{\operatorname{F}}

\newcommand{\Tb}{\mathbf{T}}

\newcommand{\Ele}{L_{2}}

\newcommand{\Elk}{\mathfrak{L}_{2}}

\newcommand{\innerp}[2]{\left\langle#1,#2\right\rangle}

\newcommand{\hsp}{{\hspace{-1pt}}}

\newcommand{\hs}{{\hspace{1pt}}}
\newcommand{\ip}[2][\cdot\hs]{\langle #1,#2\rangle}

\newcommand{\dd}{\mathrm{d}}
\newcommand{\dom}{\mathrm{dom}}

\renewcommand{\Re}{\mathcal{K}}

\newcommand{\Pn}{\operatorname{P}}
\renewcommand{\[}{\begin{equation}}
\renewcommand{\]}{\end{equation}}
\newcommand{\wegengruen}{\end{equation}}

\newcommand{\uno}{\mathbf{1}} 

\newtheorem{theorem}{Theorem}[section]
\newtheorem{corollary}[theorem]{Corollary}
\newtheorem{proposition}[theorem]{Proposition}
\newtheorem{definition}[theorem]{Definition}

\theoremstyle{definition}
\newtheorem{example}{Example}
\newtheorem{remark}{Remark}

\newcommand{\msum}[2]{\underset{{#1}}{\overset{{#2}}{\mbox{$\sum$}}}}
\newcommand{\efrac}[2]{\mbox{$\frac{#1}{#2}$}}
\numberwithin{equation}{section}

\title[Continuous frames in Krein spaces]{Continuous frames in Krein spaces}

\begin{document}
	\author[D. Carrillo]{Diego Carrillo}
\address{Departamento de Ciencias B\'asicas\\ Corporaci\'on Universitaria del Caribe\\ 
Sincelejo, Colombia.}
\email{diego.carrillo@cecar.edu.co}	
	%------author 2--------------
\author[K. Esmeral]{Kevin Esmeral}
\address{%
	Department of Mathematics\\Universidad de Caldas\\
	P.O. 170004\\
	Manizales,	Colombia.}

\email{kevin.esmeral@ucaldas.edu.co}
	%%%author 3--------------
	\author[E. Wagner]{Elmar Wagner}
	\address{Instituto de F\'isica y Matem\'aticas, Universidad Michoacana de San Nicol\'as
de Hidalgo, Edificio C-3, CU, C.P. 58040, Morelia, Michoac\'an,
M\'exico.}
	\email{elmar@ifm.umich.mx}
	\thanks{}
	\subjclass[2010]{42C15 (primary); 46B15, 46C05 (secon\-dary)}
	\keywords{Continuous frames, Krein spaces, reproducing kernels, non-regular Gram operator}
	
	\date{\today}
	
	%----------additions----------------------------------------------------------------------

\dedicatory{}
	\begin{abstract} 
	The purpose of this paper is to propose a definition of continuous frames of rank $n$ for 
		Krein  spaces and to study their basic properties. Similarly to the 
		Hilbert space case, continuous frames are characterized by the analysis,  
		the pre-frame and the frame operator, where the latter gives rise to a 
		frame decomposition theorem. The paper includes a discussion of similar, dual 
		and Parseval frames and of reproducing kernels. In addition, 
		the importance of the fundamental symmetry in the formula for the frame operator in 
		a Krein space is clarified. As prime examples, it is shown how to transfer 
		continuous frames for Hilbert spaces to Krein spaces arising from 
		a possibly non-regular Gram operator. 
	\end{abstract}
	%------------
	\maketitle
	%-------
	
	%\tableofcontents
	%%%----------------------------------------
	%\clearpage
	\section{Introduction}

A generalization of the concept of \emph{frame in a Hilbert space}  
to a family indexed by a locally compact space with a Radon measure was introduced by 
G. Kaiser \cite{Kaiser} and, independently, 
by  Ali, Antoine and Gazeau \cite{Ali-Antoine-Gazeau-1993}. 
Such frames are known as \emph{continuous frames} in distinction from
usual frames that are given by a countable spanning set.  
Continuous frames appear under different names in the literature. 
For example,  Gabardo and Han \cite{Gabardo-Han}  call them 
frames associated to measure spaces, 
and Askari-Hemmat, Dehghan, and Radjabalipour \cite{Askari} 
call them generalized frames. 
Moreover, most of the (vector) coherent states in mathematical physics 
are actually described by the continuous frames \cite{Ali-Antoine-Gazeau-2014,P}. 
For more details, see e.g. 
\cite{Ali-Antoine-Gazeau-1993,Askari,CO,Olec,Olecjensen,Chui,Foranasier,Gabardo-Han,tiro}.
   
The possibility of being over-complete makes frames more flexible 
than orthonormal bases and, for this reason, 
a powerful tool in signal processing, electrical engineering 
and several branches of mathematics 
\cite{Casazza,CasazzaLeon,CDWT,Dau,DGM,Deguang,DS,Furh,G,RNG}, 
as well as in physics 
\cite{Ali-Antoine-Gazeau-2014,Ali-Antoine-Gazeau-1993,Ali-Antoine-Gazeau1,
Ali-Antoine-Gazeau2,DGM,DS,PW}. 
On the other hand, Krein spaces arise naturally in mathematics, 
for instance in relation to signed measure spaces, 
and in mathematical physics, 
for instance in quantum field theory. 
It is therefore natural to extend  
frame theory for Hilbert spaces to Krein spaces \cite{KEFER,GMMM,K,PW}.
   
Discrete frames in a Krein space were introduced simultaneously by 
Giribet, Maestripieri, Mart\'inez Per\'ia and Massey \cite{GMMM},  
and Esmeral, Ferrer and Wagner \cite{KEFER}, albeit with slightly different definitions. 
The definition in \cite{GMMM} manifests directly the independence of a $J$-norm, 
whereas in \cite{KEFER} frames in a Krein space are essentially the same objects 
as frames for the associated Hilbert space although the property of being a frame 
also does not depend on the $J$-norm.  
Other approaches can be found in \cite{indus1,K,indus3,PW}.

The present paper introduces continuous frames of rank $n$ for 
Krein spaces by weak integrals (Definition~\ref{db})
following the Hilbert space approach developed in \cite{Ali-Antoine-Gazeau-1993}. 
As in \cite{KEFER}, continuous frames in Krein spaces are essentially the same objects 
as continuous frames in the associated Hilbert space (Theorem \ref{prop1}),  
and are characterized by the analysis,  
the pre-frame and the frame operator (see Section \ref{cf}).  
The  most important result in frame theory is the decomposition theorem 
\cite{Casazza,CasazzaLeon,Dau,DGM,Deguang,DS,G,indus2,RNG}, so it is clearly desirable 
to extend this result to the continuous case. 
This will be achieved in Theorem~\ref{Thm-trans-frame}. 
In Section~\ref{impor-fund-symm}, we discuss the importance of the fundamental 
symmetry for the invertibility of the frame operator. 
A treatment of reproducing kernels will be presented in Section~\ref{rk}. 
Unfortunately, we cannot speak about reproducing kernel Krein spaces 
unless the image of the analysis operator is ortho-complemented. 
Sections \ref{sec-dual} and \ref{sec-dual-similar} are devoted to the study of 
similar, dual and Parseval frames, where we also comment briefly on the notion of 
coherent state. In the final section, it is shown how to transfer 
continuous frames in a Hilbert space to an associated Krein space 
arising from a possibly non-regular and unbounded Gram operator $W$. 
This will be realized by extending the square root $\sqrt{|W|}$ to a 
$J$-unitary operator. The motivation for using $\sqrt{|W|}$ stems from the fact 
that a frame in the Hilbert space can never be a frame for the 
associated Krein space if the Gram operator is unbounded 
or if $0$ belongs to its spectrum.

\section{Preliminaries} \label{preliminaries}

The purpose of this section is to fix notations 
and to recall the basic elements of frame theory. 
For more details on Krein spaces, we refer the reader to \cite{Azizov} and \cite{JBognar}. 
A comprehensive introduction to frame theory can be found in \cite{Olec}. 

\subsection{Krein spaces} 
\label{KS}
Throughout this paper, $(\Re, [\cdot,\cdot])$ denotes a Krein space 
with fundamental decomposition $\Re_+\hs[\dotplus]\,\Re_-$ and 
fundamental symmetry $J$ given by 
\[ \label{J}
J (k^+ + k^-) = k^+ - k^-,\qquad  k^+ + k^-\in \Re_+\hs [\dotplus]\,\Re_-, 
\]
such that the $J$-inner product
\[ \label{ipJ}
[  h^+ \hsp +\hsp  h^-\hs, k^+\hsp  +\hsp  k^-]_J
:= [ h^+\hsp  +\hsp  h^-\hs, J (k^+\hsp  +\hsp  k^-)] 
= [  h^+ \hs, k^+ ] -[   h^-\hs,  k^-], \ \ h^\pm,k^\pm\hsp \in\hsp \Re_\pm, 
\]
turns $(\Re, [\cdot,\cdot]_J)$ into a Hilbert space. %Note that $J^2 =J$. 
The positive definite inner product $[\cdot,\cdot]_J$ defines a topology 
on $\Re$ by the $J$-norm 
\begin{equation}\label{jnorma}
\|k\|_{J}:= \sqrt{[ k\hs , k]_J} =  \sqrt{[k\hs , Jk]},\qquad k\in\Re, 
\end{equation}
and  $\Re_+\hs [\dotplus]\,\Re_-$ becomes the orthogonal sum of Hilbert spaces. 
Note that $J^2=\uno$ by \eqref{J}.

The Hilbert space  $(\mathcal{K},[\cdot,\cdot]_{J})$ is used to  study linear 
operators  acting on Krein space $(\mathcal{K},[\cdot,\cdot])$. Topological concepts 
such as continuity, closedness of operators,  spectral theory and so on refer to the 
topology induced by the\hs $J$-norm given in \eqref{jnorma}. 
Therefore, we may apply the 
same  definitions as in \emph{operator theory of Hilbert spaces}. 
The unique adjoint $\mathbf{T}^* \colon 
(\Re_2, [\cdot,\cdot]_2)\ra (\Re_1, [\cdot,\cdot]_1)$ of a 
bounded linear operator $\mathbf{T}\colon (\Re_1, [\cdot,\cdot]_1)\ra (\Re_2, [\cdot,\cdot]_2)$ 
is always taken with respect to the specified inner products, i.e., 
$$
[\mathbf{T}^* h,k  ]_1 = [h, \mathbf{T} k]_2\quad \text{for \,all}\ \ k\in\Re_1,\ \,h\in\Re_2. 
$$
Furthermore, the linear operator $\mathbf{T}^{[*]}=J_{1}\mathbf{T}^{*}J_{2}$ 
is called the $J$-\emph{adjoint} of $\mathbf{T}$, where $J_i$, $i=1,2$, denotes a (fixed) 
fundamental symmetry on $\Re_i$.  
 An operator $\mathbf{T}\in\mathcal{B}(\mathcal{K})$ is said to be \emph{self-adjoint} 
 if $\mathbf{T}=\mathbf{T}^{*}$, and $J$-\emph{self-adjoint} if
 $\mathbf{T}=\mathbf{T}^{[*]}$. A similar distinction will be made between 
 \emph{unitary} and $J$-\emph{unitary} operators. 
 By a $J$-\emph{orthogonal projection} we mean a 
 $J$-self-adjoint idempotent.

To give some examples, it follows from \eqref{J} and \eqref{ipJ} that 
$J: (\Re, [\cdot,\cdot]) \ra (\Re, [\cdot,\cdot])$ and 
$J: (\Re, [\cdot,\cdot]_J) \ra (\Re, [\cdot,\cdot]_J)$
are self-adjoint, i.e.,  $J^*=J=J^{[*]}$. 
Moreover, 
the identity operator 
\[ \label{I}
\boldsymbol{1}_J : (\Re, [\cdot,\cdot])\lra (\Re, [\cdot,\cdot]_J),\qquad 
\boldsymbol{1}_J\hs k=k, 
\]
has the adjoint $J: (\Re, [\cdot,\cdot]_J)\ra (\Re, [\cdot,\cdot])$ since 
$[ \uno_J \hs k,h]_J = [ k, h]_J =[k,J\hs  h]$
for all  $h,k\in\Re$.

A subspace $V\subset \Re$ is said to be uniformly $J$-positive 
(resp.\ uniformly $J$-negative) 
if there exists an $\varepsilon >0$ such that $[v,v]\geq \varepsilon \|v\|_J^2$ 
(resp.\ $-[v,v]\geq \varepsilon \|v\|_J^2$) for all $v\in V$. 
A self-adjoint operator $\mathbf{A}=\mathbf{A}^*$ on 
$(\Re, [\cdot,\cdot])$ is called uniformly positive, if  
$[k,\mathbf{A}k] \geq \varepsilon [ k, k]_J$ 
for a suitable constant $\varepsilon >0$ 
and all $k\in\Re$.  
Equivalently, since $[k,\mathbf{A}k]=[ k,J\mathbf{A}k]_J$, 
we have $J\mathbf{A}\geq \varepsilon$ on the Hilbert space 
$(\Re, [\cdot,\cdot]_J)$. 
As a consequence, $\mathbf{A}$ has a bounded inverse.

The fundamental projections 
\[  \label{P}
\boldsymbol{\Pn}_+:=\mbox{$\frac{1}{2}$}(\uno+J),\qquad \boldsymbol{\Pn}_-
:= \mbox{$\frac{1}{2}$}(\uno-J),\quad \uno\colon\K\rightarrow\K
\]
act on $\Re=\Re_+[\dotplus]\Re_-$ by $\boldsymbol{\Pn}_+(k^+\hsp +\hsp k^-) = k^+$, 
$\boldsymbol{\Pn}_-(k^+\hsp +\hsp k^-) = k^-$ and $\uno(k)=k$. 
Equation \eqref{P} implies immediately that  $\boldsymbol{\Pn}_\pm$ 
and  $J$ commute. Moreover, $\boldsymbol{\Pn}_+$ and $\boldsymbol{\Pn}_-$ 
are orthogonal projections, 
i.e.\ $\boldsymbol{\Pn}_{\pm}^2= \boldsymbol{\Pn}_{\pm}= \boldsymbol{\Pn}_{\pm}^*$, 
regardless of whether we 
consider $[\cdot,\cdot]$ or $[ \cdot,\cdot]_J$ on $\Re$. 

\begin{example}\label{exam-1}
	Let  $\nu$ be a signed measure on $(\mathcal{M},\mathfrak{B})$ and $\{\M_{+},\M_{-}\}$ 
	be the Hahn decomposition of $\M$ with respect to $\nu$. By the Jordan decomposition 
	theorem, there exist positive measures $\nu^{+}$ and $\nu^{-}$ on $(\M,\mathfrak{B})$ 
	such that $\nu\hsp =\hsp \nu^{+}-\nu^{-}$, \,$|\nu|\hsp =\hsp \nu^{+}+\nu^{-}$ and 
	$\nu^{+}(\M_{-})\hsp =\hsp \nu^{-}(\M_{+})\hsp =\hsp 0$.
	
	Next, let $\Elk(\M,|\nu|)$ denote the Hilbert space of 
	complex valued functions  on $\M$ that are square 
	integrable with respect to the positive measure  $|\nu|$. 
	Then, for each $f\in\Elk(\M,|\nu|)$, we may 
	write $f=\chi_{\M_{+}}f+\chi_{\M_{-}}f$, 
	where $\chi_{M}$ stands for the indicator function of $M\in \mathfrak{B}$. Thus, 
	\begin{equation}\label{decom-L2-1}
	\Elk(\M,|\nu|)=\Elk(\M_{+},\nu^{+})\,[\dotplus]\,\Elk(\M_{-},\nu^{-}).
	\end{equation}
	On $\Elk(\M,|\nu|)$,  consider the indefinite inner product 
	\begin{equation}\label{indef-metr-exam-1} 
	[\cdot,\cdot]_{\Elk}\colon\Elk(\M,|\nu|)\times\Elk(\M,|\nu|)\To\Complex, \quad 
	[f,g]_{\Elk}:=\int_{\M} \overline{f(x)} g(x)\dd\nu(x).
	\end{equation}
As easily seen,  $(\Elk(\M_{+},\nu^{+}),[\cdot,\cdot]_{\Elk})$ and 
$(\Elk(\M_{-},\nu^{-}),-[\cdot,\cdot]_{\Elk})$ are Hilbert spaces. Therefore, 
$(\Elk(\M,|\nu|),[\cdot,\cdot]_{\Elk})$ is a Krein space, 
which we abbreviate to $\Elk(\M,\nu)$. 
The projections $\mathbf{P}_{\pm}\colon\Elk(\M,\nu)\rightarrow \Elk(\M_{\pm},\nu_{\pm})$ 
are given by $\mathbf{P}_{\pm}f=\chi_{\M_{\pm}}f$ and hence the fundamental symmetry associated 
to the decomposition \eqref{decom-L2-1}, denoted by $\je$, can be expressed in terms of 
the multiplication operator with 
the Radon-Nikodym derivative $j_{\Elk}:=\chi_{\M_{+}}-\chi_{\M_{-}}$ of $\nu$ 
with respect to 
$|\nu|$, i.e., 
	\begin{equation}\label{J-L2-1}
	\je f= j_{\Elk}f= \chi_{\M_{+}}f-\chi_{\M_{-}}f,\quad f\in\Elk(\M,\nu).
	\end{equation}
	Observe that, for all $f,g\in\Elk(\M,\nu)$,  
	\begin{equation}\label{J-prod-L2}
	[ f,g]_{\je}=\int_{\M}\overline{f(x)}\hs  \je g(x)\hs \dd\nu(x)
	=\int_{\M}\overline{f(x)}\hs g(x)\hs \dd|\nu|(x)=:\innerp{f}{g}_{\Elk(\M,|\nu|)}.
	\end{equation}
	\end{example}

Let $V$ be a closed subspace of a Krein space $(\mathcal{K},[\cdot,\cdot])$. The subspace
	\begin{equation}\label{Vperp}
	V^{\perp}=\left\{x\in\mathcal{K}: [x,y]=0\,\,\text{for all $y\in\,V$}\right\}
	\end{equation}
	is  called the \emph{orthogonal complement of $V$ with respect to $[\cdot,\cdot]$}.
 A closed subspace $V$ of a Krein space $(\mathcal{K},[\cdot,\cdot])$  such that 
 $V\cap V^{\perp}=\{0\}$ and $V+ V^{\perp}=\mathcal{K}$  
	is said to be \emph{ortho-complemented}.

\begin{proposition}
\cite[Theorem 7.16]{Azizov} \label{projectivelycomplete} 
	Let $(\mathcal{K},[\cdot,\cdot])$ be a Krein space 
	and let $V$ be a closed subspace of \,$\mathcal{K}$.  
	The following statements are equivalent:
	\begin{enumerate}[i)]
		\item $V$ is ortho-complemented.
		\item $(V,[\cdot,\cdot]_{V})$ is a Krein space, where 
		$[\cdot,\cdot]_{V}:=\left.[\cdot,\cdot]\right|_{V\times V}$.
		\item Any vector in $\mathcal{K}$ admits at least 
		one $J$-orthogonal projection onto  $V$.
	\end{enumerate}
\end{proposition}

\section{Continuous frames in Krein spaces}\label{Cont-frames}

\subsection{Weak Integral in Krein spaces}

\begin{definition} \label{weakm} 
	Let $(\K_{1},[\cdot,\cdot]_{1})$, $(\K_{2},[\cdot,\cdot]_{2})$ be Krein spaces and 
	$(\M, \mathfrak{B})$  be a measureable space. 
	A function $\F\colon\M\to\mathcal{B}(\K_{1},\K_{2})$ 
	is said to be \emph{weakly measurable} in $\mathcal{B}(\K_{1},\K_{2})$ 
	if the map $\M\ni x\hs\mapsto\hs [k,\F(x)h]_{2}\in \Complex$ \,is measurable 
	for all $h\in\K_{1}$ and $k\in\K_{2}$.
\end{definition}

\begin{remark}\label{remark-weak-measu}
	For a Krein space $(\K,[\cdot,\cdot])$,  
	the Riesz representation theorem establishes an isomorphism  
	$(\K,[\cdot,\cdot]_J)\ni h \mapsto\varphi_{h}\in \K':=\mathcal{B}(\K,\mathbb{C})$, 
	where the linear functional $\varphi_{h}$ is defined by   
	$\varphi_h(\cdot)=[h,\cdot]_{J}=[h,J(\cdot)]$.  
	Thus, given a map $\mathbf{T}\colon\M\to\K$, we may consider 
	$\mathbf{U}\colon\M\to\mathcal{B}(\K,\mathbb{C})$, 
	\,$\mathbf{U}(x):=\varphi_{J\mathbf{T}(x)}$, so that 
	$\langle \lambda,\mathbf{U}(x)k\rangle_{\mathbb{C}}
	=\overline{\lambda}\hs [\mathbf{T}(x),k]$ for all $k\in \K$ and $\lambda\in\mathbb{C}$.  
	As a consequence, if $\mathbf{U}$ is weakly measurable in $\K'$, then the complex function 
	$\M\ni x\mapsto [\mathbf{T}(x),k]$ is measurable for each $k\in\K$.  
	In this case, we say that $\mathbf{T}\colon\M\to\K$ is \emph{weakly measurable} in $\K$.
\end{remark}

\begin{definition} \label{wi}
	Let $(\K_{1},[\cdot,\cdot]_{1})$, $(\K_{2},[\cdot,\cdot]_{2})$ be Krein spaces, 
	$(\M,\mathfrak{B},\nu)$ a measure space and $\F\colon\M\to\mathcal{B}(\K_{1},\K_{2})$ 
	a weakly measurable function. Then  we say that  $\F$ is \emph{weakly integrable} in 
	$\mathcal{B}(\K_{1},\K_{2})$ if there exists an 
	$\boldsymbol{A}\in\mathcal{B}(\K_{1},\K_{2})$  such that for all $h\in\K_{2}$ 
	and $k\in\K_{1}$,
	\begin{equation}\label{eg}
	\int_{\M}[h,\F(x)k]_{2}\hs \dd\nu(x)=[h,\boldsymbol{A} k]_{2}.
	\end{equation}
	The bounded operator $\boldsymbol{A}$ is called \emph{weak integral} of  \,$\F$ and 
	will be written as 
	$$\int_{\M}\F(x)\hs \dd\nu(x):=\boldsymbol{A}.$$	
\end{definition}

\begin{remark}
	Let $(\K,[\cdot,\cdot])$, $(\M,\mathfrak{B},\nu)$ and $\mathbf{T}$ be as in 
	Remark  \ref{remark-weak-measu}.  Then, 
	given any $h\in\K$ such that
	$\displaystyle\int_{\M}\langle \lambda,\mathbf{U}(x)k\rangle_{\mathbb{C}}\,\dd\nu(x)
	=\langle \lambda,\varphi_{Jh}(k)\rangle_{\mathbb{C}}$ for all $k\in \K$ and $\lambda \in\mathbb{C}$, 
	it follows that 
	$$
    \bar \lambda\!\int_{\M}[\mathbf{T}(x), k]\,\dd\nu(x)
	\,=\int_{\M}\langle \lambda,\mathbf{U}(x)\hs k\rangle_{\mathbb{C}}\,\dd\nu(x)
	=\bar \lambda\, [h,k] \,.
	$$
	Motivated by this observation, we define the weak integral of $\mathbf{T} : \M \ra \K$ in terms of 
	the weak integral of $U: \M \ra \mathcal{B}(\K,\mathbb{C})$, i.e.,  
	$$
	\int_{\M}\mathbf{T}(x)\,\dd\nu(x):=h \ \ \ \text{if and only if}\ 
	\int_{\M}[\mathbf{T}(x),k]\,\dd\nu(x)=[h,k]\ \ \ \text{for all $k\in\K$}. 
	$$
\end{remark}
The next proposition is an immediate consequence of Definition \ref{wi}. 
Note that the inte\-gra\-bility of the complex functions 
$\M\ni x\mapsto [h, F(x)k]_2\in\Complex$, \,$k\in \K_1$, \,$h\in \K_2$, implies 
the inte\-gra\-bility of $\M\ni x\mapsto [h, F(x) \mathbf{T}k_0]_2\in\Complex$ 
\,for all \,$\mathbf{T}\in\mathcal{B}(\mathcal{K}_{0},\mathcal{K}_{1})$ and $k_0\in \K_0$, 
and of $\M\ni x\mapsto [h_0, \mathbf{S}F(x)k]_3=[ \mathbf{S}^*h_0,F(x)k]_2\in\Complex$ 
\,for all \,$\mathbf{S}\in\mathcal{B}(\mathcal{K}_{2},\mathcal{K}_{3})$ and $h_0\in \K_3$, where 
$(\mathcal{K}_{i},[\cdot,\cdot]_{i})$, $i\hsp=\hsp 0,...\hs ,3$, denote Krein spaces. 
\begin{proposition}\label{pf}
Let $(\M,\mathfrak{B},\nu)$ be a measure space, 
$(\mathcal{K}_{i},[\cdot,\cdot]_{i})$, $i\hsp=\hsp 0,...\hs ,3$, 
Krein spaces,  $\F\colon \mathcal{M}\rightarrow\mathcal{B}(\mathcal{K}_{1},\mathcal{K}_{2})$ 
a weakly measurable map, and $\mathbf{A}\in\mathcal{B}(\mathcal{K}_{1},\mathcal{K}_{2})$ 
a linear operator such that $\displaystyle\int_{\mathcal{M}}\F(x)\, \dd\nu(x)=\mathbf{A}$. Then
\begin{enumerate}[i)]
\item $\displaystyle\int_{\mathcal{M}}\F(x)\hs k\,\dd\nu(x)=\mathbf{A}\hs k$ 
\;for all \,$k\in\mathcal{K}_{1}$. \\
\item  \label{pfii}
For all linear operators $\mathbf{T}\in\mathcal{B}(\mathcal{K}_{0},\mathcal{K}_{1})$ 
and $\mathbf{S}\in\mathcal{B}(\mathcal{K}_{2},\mathcal{K}_{3})$,     
$$\int_{\mathcal{M}}\F(x)\hs \mathbf{T}\,\dd\nu(x)=\mathbf{A}\mathbf{T}\quad 
\text{and}\quad\int_{\mathcal{M}}\mathbf{S}\F(x)\,\dd\nu(x)=\mathbf{S}\mathbf{A}.$$
\item $\displaystyle\int_{\mathcal{M}}\F(x)^* \dd\nu(x)=\mathbf{A}^{\hspace{-2pt}*}.$
\end{enumerate}
\end{proposition}

\subsection{Continuous frames of rank n}  \label{cf} 
From now on, we will work with a  signed measure space $(\M,\mathfrak{B},\nu)$ 
as in Example \ref{exam-1}. Our paper is based on the following fundamental definition. 
\begin{definition}\label{db} 
A  family of linearly independent vectors 
$\{\eta_{x}^{1},\eta_{x}^{2},\ldots,\eta_{x}^{n}\}_{x\in\M}$ in a 
Krein space $(\K,[\cdot,\cdot])$ is said to be a  \textit{continuous frame of rank} $n\in\N$ 
with respect to  $(\M,\mathfrak{B},\nu)$ if the functions  $\eta^{i}:\M\rightarrow\K$ given by 
$\eta^{i}(x):=\eta_{x}^{i}$ are weakly measurable for each $i=1,2,\ldots,n$ and if there exist   
positive constants $0<a\leq b$, called \emph{frame bounds}, such that 
	\begin{equation}\label{ek}
a\hs\Vert k\Vert_J^2\leq
\sum_{i=1}^n\int_{\M}|[\eta_x^i,k]|^2\hs\dd|\nu|(x)\leq b\hs \Vert k\Vert_J^2 
	\quad \text{for all }\,k\in\K.
	\end{equation}
In case $a=b$, the frame is said to be \emph{tight}, and it is called a 
\emph{Parseval frame}, if $a=b=1$. 
\end{definition}

\begin{remark}
    Clearly, the $J$-norm of a Krein space depends on the fundamental decomposition. 
    At the first glance, it seems that our definition of continuous frames of rank $n\in\N$ 
    also depends on the chosen fundamental decomposition. However, 
    two $J$-norms corresponding to different fundamental decompositions of a Krein space 
    are equivalent, see e.g. \cite[\textsection 7, Theorem 7.19]{Azizov}. 
    Therefore, only the frame bounds may change, but the property of being a continuous frame 
    in a Krein space is independent from the fundamental decomposition. 
\end{remark}

Next we state the analogue of Theorem 3.3 in \cite{KEFER} for continuous frames of rank $n$. 
This result shows that continuous frames for a Krein space are essentially the same objects 
as continuous frames for the associated Hilbert space. The proof differs from that 
in \cite{KEFER} only in notation and will be omitted. More information on the proof 
can be found in the remarks below Proposition~\ref{pb}. 
\begin{theorem}\label{prop1}
Let $n\in\N$ and let $(\K,[\cdot,\cdot])$ be a Krein space with fundamental symmetry $J$.  
Given a family of vectors  
$\{\eta_{x}^{1},\eta_{x}^{2},\ldots,\eta_{x}^{n}\}_{x\in\M}\subset\K$,  
the following statements are equivalent:
\begin{enumerate}[i)]
		\item  $\{\eta_{x}^{1},\eta_{x}^{2},...\hs,\eta_{x}^{n}\}_{x\in\M}$ 
		is a continuous frame of rank $n$ for the Krein space $(\K, [\cdot,\cdot])$ 
		with frame bounds $a\leq b$.
		
		\item  $\{\eta_{x}^{1},\eta_{x}^{2},...\hs,\eta_{x}^{n}\}_{x\in\M}$ 
		is a continuous frame of rank $n$ for the Hilbert space $(\K, [\cdot,\cdot]_{J})$ 
		with frame bounds $a\leq b$.
		
		\item  $\{J\eta_x^1,J\eta_x^2,...\hs,J\eta_x^n\}_{x\in\M}$ 
		is a continuous frame of rank $n$ for the Krein space $(\K, [\cdot,\cdot])$  
		with frame bounds $a\leq b$. %\label{iii}
		
		\item  $\{J\eta_x^1,J\eta_x^2,...\hs,J\eta_x^n\}_{x\in\M}$ 
		is a continuous frame of rank $n$ for the Hilbert space 
		$(\K, [\cdot,\cdot]_{J})$ with frame bounds $a\leq b$. %\label{iv}
\end{enumerate}
\end{theorem}

\begin{example}[Wavelets]
	Let $(L_{2}(\Real),\innerp{\cdot}{\cdot})$ be the Hilbert space consisting  
	of  all measurable functions that are square integrable with respect to the 
	Lebesgue measure of $\Real$, 
	where $\innerp{\cdot}{\cdot}$ denotes the usual inner product given by
	$\displaystyle\innerp{f}{g}=\int_{\Real}\overline{f(t)}\hs g(t)\hs\dd t$. 
	On $L_{2}(\Real)$, we consider the following sesquilinear form 
	and fundamental symmetry:   
	$$
	[f,g]:=\int_{\Real}\overline{f(t)}g(-t)\hs\dd t, \qquad (Jg)(t):=g(-t), \qquad 
	f,g\in L_{2}(\Real). 
	$$
	Note that $\innerp{f}{g}=[ f,g ]_{J}$. 
	Hence $(L_{2}(\Real),[\cdot,\cdot])$ defines a Krein space with fundamental symmetry $J$ 
	such that $(L_{2}(\Real),\innerp{\cdot}{\cdot})$ is the associated Hilbert space. 
	Let $\psi\in\,L_{2}(\Real)$ be an admisible wavelet, i.e., a function in $L_{2}(\Real)$ 
	satisfying $C_{\psi}:=\int_{\Real}\frac{|\hat{\psi}(s)|^ {2}}{|s|}\hs\dd s <\infty$, 
	where $\hat{\psi}$ denotes the Fourier transform of $\psi$. 
	For $(a,b)\in \Real\times\Real\!\setminus\!\{0\}$, set 
	$$
	\psi_{a,b}(t):=\dfrac{1}{|b|^{1/2}}\psi\left(\dfrac{t-a}{b}\right),\quad t\in\Real.
	$$
It is well known that $\{\psi_{a,b}\}_{(a,b)\in\Real\times\Real\setminus\{0\}}$ defines  
a tight continuous frame of rank $1$ with frame bound $C_{\psi}$ for the Hilbert space 
$(L_{2}(\Real),\innerp{\cdot}{\cdot})$ and with respect to 
$\left(\Real\times\Real\!\setminus\!\{0\},\mathfrak{B},\frac{\dd a\hs \dd  b}{b^{2}}\right)$, 
see \cite[Corollary 11.1.2]{Olec}.  From Proposition \ref{prop1}, we conclude that 
$\{\psi_{a,b}\}_{(a,b)\in\Real\times\Real\setminus\{0\}}$  yields a tight continuous frame 
for the Krein space $(L_{2}(\Real),[\cdot,\cdot])$, with the same frame bound $C_{\psi}$ 
and with respect to the same measure space.
\end{example}

%----------frame transform---------------- 
Our next aim is to define the analogue of the frame operator for continuous frames. 
Let $\{\eta_{x}^{1},\eta_{x}^{2},\ldots,\eta_{x}^{n}\}_{x\in\M}$ be a continuous frame 
of rank $n\in\N$ for the Krein space $(\K,[\cdot,\cdot])$ with frame bounds $0<a\leq b$ 
and with respect to the measure space $(\M,\mathfrak{B},\nu)$. 
Combining Definition~\ref{weakm}, Remark \ref{remark-weak-measu} and Definition~\ref{db} 
shows that the function of rank 1 operators 
$\M \ni x\mapsto |\eta_x^{i}][\eta_x^{i}| \in \mathcal{B}(\K)$ is weakly measurable, where 
$|\eta_x^{i}][\eta_x^{i}|(k):=[\eta_x^{i},k]\hs \eta_{x}^{i}$ for $k\in\K$ and 
$i=1,\ldots,n$. 
Consider the 
the sesquilinear form $\Psi\colon\K\times\K\To \Complex$ given by 
$$
\Psi(k_{1},k_{2})=\sum_{i=1}^{n}\int_{\M}
[k_{1},\eta_{x}^{i}]\hs[\eta_{x}^{i},k_{2}]\hs\dd|\nu|(x),\quad k_{1},k_{2}\in\K.
$$
By the Cauchy-Schwarz inequality, we have  
$\left|\Psi(k_{1},k_{2})\right|\leq b\hs \|k_{1}\|_{J}\hs\|k_{2}\|_{J}$ 
for all $k_{1},k_{2}\in\K$,  
so $\|\Psi\|\leq b$. Thus, the Riesz representation theorem ensures that there exists a unique 
bounded linear operator $\mathbf{S}\colon\K\To\K$ such that  $\|\mathbf{S}\|=\|\Psi\|$ 
and $\Psi(k_{1},k_{2})=[k_{1},\mathbf{S}k_{2}]$. Since 
$$
[k_1,\mathbf{S}k_2] 
=\sum_{i=1}^{n}\int_{\M}
[k_{1},\eta_{x}^{i}]\hs[\eta_{x}^{i},k_{2}]\hs\dd|\nu|(x) 
=\sum_{i=1}^{n}\int_{\M}[k_1,j_{\Elk}\hsp (x)\hs \eta_{x}^{i}]\,[\eta_{x}^{i},k_2]\,\dd\nu(x), 
$$
where $j_{\Elk}:=\chi_{\M_{+}}\!-\chi_{\M_{-}}$ denotes the Radon-Nikodym derivative 
of $\nu$ with respect to $|\nu|$, 
it follows from Definition \ref{wi} that 
\begin{equation}\label{transformada-frame}
\mathbf{S}=\sum_{i=1}^{n}\int_{\M}|\eta_{x}^{i}]\hs [\eta_{x}^{i}|\,\dd|\nu|(x)
=\sum_{i=1}^{n}\int_{\M}j_{\Elk}\hsp (x)|\eta_{x}^{i}]\hs [\eta_{x}^{i}|\,\dd\nu(x)
\end{equation} 
is the weak integral of the weakly measurable function 
$\M \ni x\mapsto |\eta_x^{i}][\eta_x^{i}| \in \mathcal{B}(\K)$ with respect to $|\nu|$,  
and of $\M \ni x\mapsto j_{\Elk}\hsp (x)\hs |\eta_x^{i}] [\eta_x^{i}| \in \mathcal{B}(\K)$ 
with respect to $\nu$. 

We call the weak integral $\mathbf{S}$ in \eqref{transformada-frame} 
the \emph{frame operator} and $\mathbf{S}k$ the \emph{frame transform} of $k\in\K$.
From \eqref{ek} and \eqref{transformada-frame}, it follows that 
\begin{equation}\label{desi:trans-Frame}
a\,[k,k]_{J}\,\leq\,[ k,\mathbf{S}k]\,\leq\,b\,[k,k]_{J}\quad \text{for \,all \,$k\in\K$}.
\end{equation}
As a consequence, $\mathbf{S}$ is uniformly positive and therefore invertible. 
Moreover, the inverse operator $\mathbf{S}^{-1}$, 
regarded as an operator on the Hilbert space $(\Re, [\cdot,\cdot]_J)$, satisfies
\[ \label{Sinv} 
0<b^{-1}\uno \leq \mathbf{S}^{-1}J\leq a^{-1}\uno \quad \text{and} \quad 
0<b^{-1}\uno \leq J\mathbf{S}^{-1}\leq a^{-1}\uno, 
\] 
where the second relation follows from the first by applying $J^2=\uno$. 

To describe the frame operator in terms of an analogue of the 
so-called analysis operator, 
let $\bigoplus_{i=1}^n\Elk(\M,\nu)$ be the $n$-fold direct orthogonal sum of the 
Krein space $\Elk(\M,\nu)$ from Example \ref{exam-1}. 
For $f\in \bigoplus_{i=1}^n\Elk(\M,\nu)$, we write $f=(f_i)_{i=1}^{n}$  and 
$f(x)=\big(f_i(x)\big)_{i=1}^{n}$, where $f_i\in\Elk(\M,\nu)$. 
Consider the inner product $[\cdot,\cdot]^{\oplus}$ in 
\,$\bigoplus_{i=1}^n\!\Elk(\M,\nu)$ given by  	
$$
[(g_{i})_{i=1}^{n},(f_i)_{i=1}^{n}]^{\oplus}
:=\sum_{i=1}^n\int_{\M}\overline{g_i(x)}\hs f_i(x)\hs\dd\nu(x)\hs .
$$
Then $\left(\bigoplus_{i=1}^n\!\Elk(\M,\nu)\hs,\hs[\cdot,\cdot]^{\oplus}\right)$ 
yields a Krein space with fundamental decomposition $$
\bigoplus_{i=1}^n\Elk(\M,\nu)
=\bigoplus_{i=1}^n\Elk(\M_{+},\nu_+)\,[\dotplus]\,\bigoplus_{i=1}^n\Elk(\M_{-},\nu_-)
$$
and fundamental symmetry 
$J_{\hsp\Elk}^\oplus\big((f_{i})_{i=1}^{n}\big):=(J_{\Elk}f_{i})_{i=1}^{n}$, 
where $J_{\Elk}$ denotes the fundamental symmetry of the Krein space 
$(\Elk(\M,\nu),[\cdot,\cdot])$ as described in \eqref{J-L2-1}.
The associated Hilbert space 
$\Big(\bigoplus_{i=1}^n\!\Elk(\M,\nu)\hs,\hs
[\cdot,\cdot]^{\oplus}_{J_{\hsp\Elk}^\oplus}\Big)$
has the inner product 
\[  \label{ipL2}
[ (g_{i})_{i=1}^{n},(f_{i})_{i=1}^{n}]_{J_{\hsp\Elk}^\oplus}^\oplus
=[(g_{i})_{i=1}^{n}, J_{\hsp\Elk}^\oplus \big((f_{i})_{i=1}^{n}\big)]^{\oplus}
=\sum_{j=1}^{n}\innerp{g_{j}}{f_{j}}_{\Elk(\M,|\nu|)}.
\]
Define 
\[\label{oper-T[*]}
\mathbf{T}\colon \K\rightarrow \bigoplus_{i=1}^{n}\Elk(\M,\nu), \quad 
(\mathbf{T}k)(x):=\left([\eta_{x}^{i},k]\right)_{i=1}^{n},\quad k\in\K.
\]
It follows immediately from \eqref{ek} that $\mathbf{T}$ is well defined and bounded. 
In analogy to the Hilbert space case 
\cite{Ali-Antoine-Gazeau-2014,Ali-Antoine-Gazeau-1993,Azhini-Beheshti}, 
we refer to $\mathbf{T}$ as the \emph{analysis operator}, to $\mathbf{T}k$ as the 
 \emph{the analysis transform of} $k\in \K$, and to 
$\widehat{k}_{i}(x):=[\eta_{x}^{i},k]$ as the \emph{$i$-th analysis transform of} $k$. 
Straightforward calculations show that 
its adjoint operator, the so-called \emph{pre-frame operator}, is given by the weak integral 
\[\label{pre-frame-op} 
\mathbf{T}^{*}(f_i)_{i=1}^{n}
=\sum_{i=1}^n\int_{\M}f_i(x)\hs\eta_{x}^i\hs\dd\nu(x),\qquad f_i\in\Elk(\M,\nu).
\]
With the bounded operators 
$\mathrm{T}_i^{*}\colon\Elk(\M,\nu)\rightarrow\K$, 
\,\hs$\mathrm{T}_i^{*} f_i:=\int_{\M}f_i(x)\hs\eta_x^i\hs\dd\nu(x)$, 
\,Equation~\eqref{pre-frame-op} becomes 
$\mathbf{T}^{*}(f_{i})_{i=1}^{n}=\msum{i=1}{n}\hs \mathrm{T}_{i}^{*}f_{i}$. 
Furthermore, $\mathbf{T}^{*} J_{\hsp\Elk}^\oplus \mathbf{T}k 
=\msum{i=1}{n}\int_{\M}j_{\Elk}\hsp (x)\hs [\eta_{x}^{i},k]\hs 
\eta_{x}^{i}\,\dd\nu(x)$
for all $k\in \K$, 
so that
\[ \label{Sad}
\mathbf{S} = \mathbf{T}^{*} J_{\hsp\Elk}^\oplus \mathbf{T} 
\]
by \eqref{transformada-frame}.  

The most important result in frame theory is the frame decomposition theorem. 
For continuous frames, it states that each vector in the Krein (Hilbert) space 
can be represented by a weak integral obtained from the continuous frame 
and the inverse of the 
frame operator. An equivalent way of stating the frame decomposition theorem is that 
this weak integral yields the identity operator. 

\begin{theorem}[{\bf Frame decomposition}]\label{Thm-trans-frame}
	Let $(\K,[\cdot,\cdot])$ be a Krein space with fundamental symmmetry $J$ 
	and let $\{\eta_{x}^{1},\eta_{x}^{2},\ldots,\eta_{x}^{n}\}_{x\in\M}$ 
	be a continuous frame of rank $n\in\N$ for $\K$ 
	with respect to $(\M,\mathfrak{B},\nu)$. 
	Then 
\begin{align}
\uno_{\K}&=     \label{sinf2}
\sum_{i=1}^{n}\int_{\M}|\eta_{x}^{i}]\,[\mathbf{S}^{-1}\eta_{x}^{i}|\,\dd|\nu|(x)
=\sum_{i=1}^{n}\int_{\M}|\mathbf{S}^{-1}\eta_{x}^{i}]\,[\eta_{x}^{i}|\,\dd|\nu|(x).
\end{align}
\end{theorem}
\begin{proof}
It has been observed in \eqref{Sinv} that $\mathbf{S}$ has a bounded inverse. 
Clearly, by \eqref{Sad}, we have $\mathbf{S}^*=\mathbf{S}$ and thus 
$\mathbf{S}^{-1*}=\mathbf{S}^{-1}$. 
From \eqref{transformada-frame} and Proposition \ref{pf}.\ref{pfii}), it follows that 
$$
k= \mathbf{S}^{-1} \mathbf{S} k 
=\sum_{i=1}^{n}\int_{\M} [\eta_{x}^{i},k]\hs\,\mathbf{S}^{-1}\eta_{x}^{i} \,\dd|\nu|(x)
\ \,\text{ and } \ \,
k=  \mathbf{S} \mathbf{S}^{-1} k 
=\sum_{i=1}^{n}\int_{\M} [\mathbf{S}^{-1}\eta_{x}^{i},k]\hs\,\eta_{x}^{i} \,\dd|\nu|(x)
$$
for all $k\in \K$. Now the definition of the weak integral gives the result. 
\end{proof}
   
By Theorem \ref{prop1}, each continuous  frame 
$\{\eta_{x}^{1},\eta_{x}^{2},\ldots,\eta_{x}^{n}\}_{x\in\M}$ 
of rank $n$ for  $(\K, [\cdot,\cdot])$ 
with respect to $(\M,\mathfrak{B},\nu)$ gives rise to 
three other frames with slightly different frame operators. 
In the following, we will relate these frame operators to $\mathbf{S}$ 
presented in \eqref{transformada-frame}. 
First, consider the frame $\{J\eta_x^1,J\eta_x^2,\ldots,J\eta_x^n\}_{x\in\M}$ 
for $(\K, [\cdot,\cdot])$. 
Denoting the corresponding frame operator  by $\mathbf{S}_{0}$, 
we get from \eqref{transformada-frame}
\begin{equation}\label{S0}
\mathbf{S}_{0}\hs k = \sum_{i=1}^{n}\int_{\M} [J\eta_{x}^{i},k]\hs J\eta_{x}^{i}\hs\dd|\nu|(x)
= \sum_{i=1}^{n}\int_{\M} [\eta_{x}^ {i},Jk]\hs J\eta_{x}^{i}\hs\dd|\nu|(x),\quad k\in\K.
\end{equation}
Comparing \eqref{S0} with \eqref{transformada-frame} 
and applying Proposition \ref{pf}.\ref{pfii}) shows that the two frame operators 
are related by 
$
\mathbf{S}_{0} = J\hs \mathbf{S}\hs J. 
$
Next, let $\mathbf{S}_{1}$ be the frame operator of the frame 
$\{\eta_{x}^{1},\eta_{x}^{2},\ldots,\eta_{x}^{n}\}_{x\in\M}$ 
for the Hilbert space $(\K, [\cdot,\cdot]_{J})$. 
Then, by \eqref{transformada-frame}, we have 

\begin{equation}\label{S1}
 \mathbf{S}_{1}\hs k 
 = \sum_{i=1}^{n}\int_{\M} [\eta_{x}^ {i},k]_J\hs \eta_{x}^{i}\hs\dd|\nu|(x) 
 = \sum_{i=1}^{n}\int_{\M} [\eta_{x}^{i},J k]\hs \eta_{x}^{i}\dd|\nu|(x), \quad k\in\K,  
\end{equation}
and thus 
$
\mathbf{S}_{1}= \mathbf{S}J. 
$
Finally, with $\mathbf{S}_{2}$ denoting the frame operator of the continuous frame 
$\{J\eta_x^1,J\eta_x^2,\ldots,J\eta_x^n\}_{x\in\M}$ 
for the Hilbert space $(\K, [\cdot,\cdot]_{J})$, we get from \eqref{transformada-frame}
\begin{equation}\label{S2}
\mathbf{S}_{2}k 
= \sum_{i=1}^{n}\int_{\M} [J\eta_{x}^{i},k]_{J}\hs J \eta_{x}^{i} \hs\dd|\nu|(x)
= \sum_{i=1}^{n}\int_{\M} [\eta_{x}^{i},k]\hs J \eta_{x}^{i}\hs\dd|\nu|(x), \quad k\in\K, 
\end{equation}
so that 
$
\mathbf{S}_{2} = J \hs \mathbf{S}. 
$

The next two propositions show that 
continuous frames can be characterized by properties of the analysis and  
pre-frame operator. By Theorem \ref{prop1}, the proofs are literally the same 
as in the Hilbert space case, see e.g.\ \cite{Kaiser}.

\begin{proposition}\label{prop:injectivity-analysis-operator} 
Let $(\M,\mathfrak{B},\nu)$ be a measure space and let 
$(\K,[\cdot,\cdot])$ be a Krein space. Given a collection of weakly measurable functions  
$\M\ni x\mapsto \eta^j_x\in\K$, \,$j=1,\ldots,n$, the family 
$\{\eta_{x}^{1},\eta_{x}^{2},\ldots,\eta_{x}^{n}\}_{x\in\M}$ is a continuous frame 
of rank $n\in\N$ for $\K$ with respect to $(\M,\mathfrak{B},\nu)$  if and only if the 
analysis operator $\mathbf{T}$ defined in \eqref{oper-T[*]} is injective.
\end{proposition}

\begin{proposition}\label{prop:sur-pre-frame-operator} 
Let $(\M,\mathfrak{B},\nu)$ be a measure space and let 
$(\K,[\cdot,\cdot])$ be a Krein space. Given a collection of weakly measurable functions  
$\M\ni x\mapsto \eta^j_x\in\K$, \,$j=1,\ldots,n$, the family 
$\{\eta_{x}^{1},\eta_{x}^{2},\ldots,\eta_{x}^{n}\}_{x\in\M}$ is a continuous frame 
of rank $n\in\N$ for $\K$ with respect to $(\M,\mathfrak{B},\nu)$  if and only if the 
pre-frame operator $\mathbf{T}^{*}$ given in \eqref{pre-frame-op} is surjective.
\end{proposition}

In practice, it is not necessary to verify the frame condition \eqref{ek} for all $k\in\K$, 
it suffices to show \eqref{ek} on a dense subspace. This will be proven in the next proposition. 

\begin{proposition}\label{teo-impo-no-nonregular}
Let $(\M,\mathfrak{B},\nu)$ be a measure space and let 
$\{\eta_{x}^{1},\eta_{x}^{2},\ldots,\eta_{x}^{n}\}_{x\in\M}$ be a family of vectors in a Krein 
space $(\K,[\cdot,\cdot])$ such that the maps 
$\M\ni x\mapsto\eta_{x}^i\in\K$, \,$i\hsp =\hsp 1,\ldots,n$, 
are weakly measurable. If the frame condition \eqref{ek} is satisfied for all elements 
in a dense subset $\G$ of $\K$, 
then $\{\eta_{x}^{1},\eta_{x}^{2},\ldots,\eta_{x}^{n}\}_{x\in\M}$ 
is a continuous frame of rank $n\in\N$ for $\K$ with respect to $(\M,\mathfrak{B},\nu)$.
\end{proposition} 

\begin{proof} 
Assume that  \eqref{ek} is satisfied for all $k\in\G$ and for 
fixed positive constants $a\leq b$.  
Using \eqref{ipL2} and \eqref{oper-T[*]}, we can write \eqref{ek} as 
\[ \label{Tfb} 
a\hs \Vert k\Vert_J^2\ \leq\ \Vert \mathbf{T}k\Vert_{J_{\hsp\Elk}^\oplus}^2
\ \leq\  b\hs \Vert k\Vert_J^2  \quad \text{for \,all }\,k\hs\in\hs \G. 
\]
It follows that $\mathbf{T}$ is bounded and thus uniformly continuous. 
By uniform continuity, \eqref{Tfb} holds for all $k\in\K$, which is equivalent to \eqref{ek}. 
\end{proof}

We close this section by showing that any continuous frame of rank $n$ can be reduced to 
a continuous frame of rank $1$ with respect to a suitable the measure space.

\begin{proposition} \label{r1}
Any continuous frame  $\{\eta_{x}^{1},\eta_{x}^{2},\ldots,\eta_{x}^{n}\}_{x\in\M}$ 
of rank $n\in\N$ for a Krein space $(\K,[\cdot,\cdot])$ with respect to 
$(\M,\mathfrak{B},\nu)$ may be viewed as a continuous frame 
$\{\sigma(y)\}_{y\in\M_{n}}$ of rank 1 for $\K$ 
with respect to some measure space $(\M_{n},\mathfrak{B}_{n},\nu_{n})$.
\end{proposition}	
\begin{proof}
Assume that $\{\eta_{x}^{1},\eta_{x}^{2},\ldots,\eta_{x}^{n}\}_{x\in\M}$ 
is a continuous frame of rank $n\in\N$ for the Krein space $(\K,[\cdot,\cdot])$ 
with respect to $(\M,\mathfrak{B},\nu)$. 
Set $\M_n:=\{1,\ldots,n\}\times \M$ and consider the counting measure $\mu_n$ 
on the power set $2^{\{1,\ldots,n\}}$. 
Let $\nu_n:= \mu_n\otimes \nu$ denote the product measure on 
product $\Sigma$-algebra $\mathfrak{B}_n:= 2^{\{1,\ldots,n\}}\otimes\mathfrak{B}$. 
Define $\sigma: \M_n \ra \K$ by $\sigma((i,x)):=\eta^i_x$. Since 
the functions $\M\ni x\mapsto \eta_x^i\in \K$, \,$i=1,\ldots,n$, 
are weakly measurable, 
it follows that $\sigma$ is weakly measurable. 
Note that $|\nu_n|=\mu_n\otimes |\nu|$. 
Now,  by Tonelli's theorem, we have 
\begin{align*}
\int_{\M_{n}}|[\sigma(i,x),k]|^2\,d|\nu_{n}|((i,x))
&=\sum_{i=1}^n\int_{\M}
|[\eta_x^i,k]|^2\,d|\nu|(x)
\end{align*}
for all $k\in \K$. From this, we conclude that $\{\sigma(y)\}_{y\in\M_{n}}$ 
is a continuous frame of rank 1 for $\K$ 
with respect to $(\M_{n},\mathfrak{B}_{n},\nu_{n})$  
admitting the same frame bounds as 
$\{\eta_{x}^{1},\eta_{x}^{2},\ldots,\eta_{x}^{n}\}_{x\in\M}$. 
\end{proof}

%------importance of fund-simetry------
%\clearpage

\subsection{The role of fundamental symmetry in the formula of frame 
operator}\label{impor-fund-symm}
Given  a continuous frame $\{\eta_{x}^{1},\eta_{x}^{2},\ldots,\eta_{x}^{n}\}_{x\in\M}$
of rank $n$ for the Krein space $(\K, [\cdot,\cdot])$, 
let $\mathbf{T}$ and $\mathbf{T}^{*}$  denote the 
analysis operator and pre-frame operator defined 
in \eqref{oper-T[*]} and \eqref{pre-frame-op}, respectively. 
In this section, we will discuss the importance of the   fundamental symmetry 
$J_{\hsp\Elk}^\oplus$ in the formula \eqref{Sad} 
of the frame operator: $\mathbf{S}=\mathbf{T}^{*}J_{\hsp\Elk}^\oplus\mathbf{T}$.  

For comparison, consider the bounded linear operator 
$\s\colon\K\rightarrow\K$ given by  $\s:=\mathbf{T}^{*}\mathbf{T}$. 
As the frame operator $\mathbf{S}$, it is obviously self-adjoint 
in the Krein space $(\K, [\cdot,\cdot])$. 
Furthermore, by \eqref{oper-T[*]} and \eqref{pre-frame-op}, 
\begin{equation}\label{transformada-Stilde}
    [ h,\s k]=\sum_{i=1}^{n}\int_{\M}[\eta_{x}^{i},k]\hs [h,\eta_{x}^{i}]\,\dd\nu(x) 
\end{equation}
for all $h,k\in \K$. Hence $\s$ coincides with the weak integral
$$
\s =\sum_{i=1}^{n}\int_{\M}|\eta_{x}^ {i}]\hs[\eta_{x}^ {i}|\hs\hs\dd\nu(x).
$$

The operators $\s$ and $\mathbf{S}$ have different properties,  
for instance, $\mathbf{S}$ is always invertible, 
but $\s$ may not be and therefore cannot be used in the frame decomposition theorem.  
We illustrate this in the following example.

\begin{example}\label{examp-import}
Let $(\mathcal{M}_0,\mathfrak{B}_0,\mu_0)$ be a finite measure space and 
$\{\eta_x^1,\eta_x^2,\ldots,\eta_x^n\}_{x\in\mathcal{M}_0}$ 
a continuous frame of rank $n\in\N$ for a Krein space $(\K,[\cdot,\cdot])$ 
with respect to $(\mathcal{M}_0,\mathfrak{B}_0,\mu_0)$ 
and with frame bounds $0<a\leq b$. Consider a 
two point set $\{p_-,p_+\}$ and the signed measure $\upsilon$ on the power set 
$2^{\{p_-,p_+\}}$ determined by $\upsilon(\{p_-\})=-1$ and $\upsilon(\{p_+\})=1$. 
Let $\mathfrak{B}:=\mathfrak{B}_0\otimes 2^{\{p_-,p_+\}}$ denote the product $\Sigma$-algebra 
and $\nu:=\mu_0\otimes\upsilon$ the product measure on $\mathfrak{B}$. 
For $(x,p_{\pm})\in \mathcal{M}:=\mathcal{M}_0\times\{p_-,p_+\}$, set 
$\eta_{(x,p_{\pm})}^i:=\eta_x^i$. Then the map   
$\mathcal{M}\ni (x,p) \mapsto \eta_{(x,p)}^i\in \K$ is weakly measurable 
for all $i=1,\ldots,n$ and, by Tonelli's Theorem, 
\begin{align*} 
\sum_{i=1}^n\int_{\mathcal{M}}|[\eta_{(x,p)}^i,k]|^2\,\dd|\nu|((x,p)) 
&= \sum_{i=1}^n \int_{\mathcal{M}_0}\left(
\int_{\{p_{\pm}\}}|[\eta_{(x,p)}^i,k]|^2\dd|\upsilon|(p)\right)\dd|\mu_0|(x) \\
 & =  2\sum_{i=1}^n\int_{\mathcal{M}_0}|[\eta_x^i,k]|^2\,\dd|\mu_0|(x).
\end{align*} 
Hence $\{\eta_{(x,p)}^1,\eta_{(x,p)}^2,\ldots\eta_{(x,p)}^n\}_{(x,p)\in\mathcal{M}}$ is a 
continuous frame of rank $n\in\N$ for the Krein space $\K$ 
with frame bounds $0<2a\leq 2b$.  Now, by Fubini's theorem, 
Equation \eqref{transformada-Stilde} yields 
\begin{align*}
[ h,\s k]&=\sum_{i=1}^{n}\int_{\M}[\eta_{(x,p)}^i,k]\hs [h,\eta_{(x,p)}^i]\,\dd\nu((x,p))\\ 
&=\sum_{i=1}^{n}\int_{\M_0}\left([\eta_{(x,p)}^i,k]\hs [h,\eta_{(x,p)}^i] 
-[\eta_{(x,p)}^i,k]\hs [h,\eta_{(x,p)}^i]\right)\,\dd\mu_0(x) \ =\ 0 
\end{align*}
for all $h,k\in\K$. 
As $(\K,[\cdot,\cdot])$ is non-degenerate, it follows that  
$\s=\mathbf{T}^*\mathbf{T}=0$. 
In particular, $\s$ is not invertible.
\end{example}

\subsection{Reproducing kernels} \label{rk}
Let $\{\eta_{x}^{1},\eta_{x}^{2},...\hs,\eta_{x}^{n}\}_{x\in\M}$ 
be a continuous frame of rank $n$ for the Krein space $(\K, [\cdot,\cdot])$ 
with frame bounds $a\leq b$. Then, by Theorem \eqref{prop1}, we may view 
$\{\eta_{x}^{1},\eta_{x}^{2},...\hs,\eta_{x}^{n}\}_{x\in\M}$ as a continuous frame 
in the Hilbert space  $(\K, [\cdot,\cdot]_J)$ with the same frame bounds. 
Considering the analysis operator $\mathbf{T}$ 
defined in \eqref{oper-T[*]} as a Hilbert space operator 
$\mathbf{T} : (\K, [\cdot,\cdot]_J) \lra\Big(\bigoplus_{i=1}^n\!\Elk(\M,\nu)\hs,\hs
[\cdot,\cdot]^{\oplus}_{J_{\hsp\Elk}^\oplus}\Big)$, it follows from \eqref{Tfb} that 
$\mathbf{T}$ is bounded, injective, and $\Rang \Tb$ is closed. In particular, 
$\Tb(\K):=\Rang \Tb$ with the inner product $[\cdot,\cdot]^{\oplus}_{J_{\hsp\Elk}^\oplus}$ 
is a Hilbert space by itself. As shown in 
\cite{Ali-Antoine-Gazeau-2014} and 
\cite{Ali-Antoine-Gazeau1}, 
$\Big(\Tb(\K),[\cdot,\cdot]^{\oplus}_{J_{\hsp\Elk}^\oplus}\Big)$ 
is actually a reproducing kernel Hilbert space. To see this, note that 
$$
\left([\eta_{x}^{i},k]\right)_{i=1}^{n} 
=\Big(\sum_{j=1}^{n}\int_{\M}\!
[\eta_{x}^{i},\mathbf{S}^{-1}\eta_{y}^{j}]\,[\eta_{y}^{j},k]\,\dd|\nu|(y) \Big)_{i=1}^{n} 
\quad \text{for \,all }\,k\in \K,\ \,x\in \M,  
$$
by \eqref{sinf2}. Setting $K_{ij}(x,y):= [\eta_{x}^{i},\mathbf{S}^{-1}\eta_{y}^{j}]$ 
and $\mathbf{K}(x,y):= \big(K_{ij}(x,y)\big)_{i,j=1}^n
\in\mathrm{Mat}_{n\times n}(\mathbb{C})$, we have for all 
$\Phi=(\phi_i)_{i=1}^n\in\Tb(\K)$ 
\[ \label{Phix}
\Phi(x)= \int_{\M}\! \mathbf{K}(x,y)\circ \Phi(y)\,\dd|\nu|(y)  
=\Big(\sum_{j=1}^{n}\int_{\M}\! K_{ij}(x,y)\, \phi_j(y)\,\dd|\nu|(y) \Big)_{i=1}^{n},  
\]
where $\circ$ denotes matrix multiplication.  Moreover, since 
\[ \label{barK} 
\overline{K_{ij}(x,y)}=[\eta_{y}^{j},\mathbf{S}^{-1}\eta_{x}^{i}]=K_{ji}(y,x), 
\]
it follows that 
$$
\| \Phi(x)\|_{\mathbb{C}^n}^2 \leq 
\Big(\sum_{i=1}^{n}\sum_{j=1}^{n}\int_{\M}\![\eta_{x}^{i},\mathbf{S}^{-1}\eta_{y}^{j}]\,
[\eta_{y}^{j},\mathbf{S}^{-1}\eta_{x}^{j}]\,\dd|\nu|(y) \Big)\, \|\Phi\|_{J_{\hsp\Elk}^\oplus}
=\Big(\sum_{i=1}^{n} [\eta_{x}^{j},\mathbf{S}^{-1}\eta_{x}^{j}] \Big)\,
\|\Phi\|_{J_{\hsp\Elk}^\oplus}. 
$$
As a consequence, the evaluation map $\Tb(\K)\ni \Phi \mapsto \phi_i(x) \in \mathbb{C}$ 
is continuous for all $x\in \M$ and $i=1,\ldots,n$. 
Moreover, Equations \eqref{oper-T[*]}, \eqref{Phix}  and \eqref{barK} show that 
this evaluation map can be given by the following inner product: 
\[ \label{phiix}
\phi_i(x) = [\Tb(\mathbf{S}^{-1}\eta_{x}^{i}) ,\Phi]^{\oplus}_{J_{\hsp\Elk}^\oplus}= 
[\big(\hs\overline{K_{ij}(x,\cdot)}\hs\big)_{j=1}^{n},
\Phi(\cdot)]^{\oplus}_{J_{\hsp\Elk}^\oplus}. 
\]

Since we are working in the Krein space setting, it would be desirable to have 
a Krein space analogue reproducing kernel Hilbert space. 
Clearly, \eqref{phiix} can be written 
\[  \label{Kphiix}
\phi_i(x) = [J_{\hsp\Elk}^\oplus\Tb(\mathbf{S}^{-1}\eta_{x}^{i}) ,\Phi]^{\oplus} = 
[\big(j_{\Elk}(\cdot)\,\overline{K_{ij}(x,\cdot)}\hs\big)_{j=1}^{n},\Phi(\cdot)]^{\oplus}, 
\]
and then \eqref{Phix} becomes 
\[  \label{KPhix}
\Phi(x)= \int_{\M}\! J_{\hsp\Elk}^\oplus\mathbf{K}(x,y)\circ \Phi(y)\,\dd \nu(y)  
=\Big(\sum_{j=1}^{n}\int_{\M}\! j_{\Elk}(y)\,K_{ij}(x,y)\, \phi_j(y)\,\dd \nu(y) 
\Big)_{i=1}^{n}. 
\]
However, although $\Tb(\K)$ is always a \emph{closed} subspace of 
$\Big(\bigoplus_{i=1}^n\!\Elk(\M,\nu)\hs,\hs
[\cdot,\cdot]^{\oplus}\Big)$, the restriction of $[\cdot,\cdot]^{\oplus}$ to 
$\Tb(\K)$ does not necessarily yield a Krein space. 
For instance, in Example~\ref{examp-import}, we have 
$[\Tb h,\Tb k]=[ h,\s k]=0$ for all $h,k\in \K$ so that 
$\Tb(\K) \cap \Tb(\K)^\bot = \Tb(\K)$. 
Therefore, by Proposition \eqref{projectivelycomplete}, we may view 
$\big(\Tb(\K),[\cdot,\cdot]^{\oplus}\big)$ as a 
\emph{reproducing kernel Krein space} only if 
any element in $\bigoplus_{i=1}^n\!\Elk(\M,\nu)$ admits at least 
one $J_{\hsp\Elk}^\oplus$-orthogonal projection onto $\Tb(\K)$, or equivalently, if
$\Tb(\K)$ is ortho-complemented. In this case, \eqref{Kphiix} and \eqref{KPhix} 
hold with the Krein space kernel function 
$$
J_{\hsp\Elk}^\oplus\mathbf{K}: \M \times \M \ra \mathrm{Mat}_{n\times n}(\mathbb{C}),\ \ 
J_{\hsp\Elk}^\oplus\mathbf{K}(x,y)=
\left(j_{\Elk}(y)\hs K_{ij}(x,y)\right)_{i,j=1}^n
\!=\left(j_{\Elk}(y)\hs [\eta_{x}^{i},\mathbf{S}^{-1}\eta_{y}^{j}] \right)_{i,j=1}^n.  
$$
Moreover, if $J_{\hsp\Elk}^\oplus : \Tb(\K) \ra \Tb(\K)$, then 
$$
K_x^i:=J_{\hsp\Elk}^\oplus \Tb(\mathbf{S}^{-1}\eta_{x}^{i}) 
= \left(j_{\Elk}(\cdot)\,\overline{K_{ij}(x,\cdot)}\right)_{j=1}^{n} \in \Tb(\K) 
$$
and 
$$
\sum_{i=1}^n\int_{\M}\left|[K_x^i,\Phi]^\oplus\right|^2\hs\dd|\nu|(x) 
= \sum_{i=1}^n \int_{\M}\Big| \sum_{j=1}^n \int_{\M} 
\! K_{ij}(x,y)\, \phi_j(y)\,
\hs\dd|\nu|(y) \Big|^2 \dd|\nu|(x) 
= \|\phi\|_{J_{\hsp\Elk}^\oplus}^2, 
$$
so $\{K_{x}^{1},K_{x}^{2},\ldots,K_{x}^{n}\}_{x\in\M}$
defines a Parseval frame for $\big(\Tb(\K),[\cdot,\cdot]^{\oplus}\big)$.

%----dual and similar-----

\subsection{Similar continuous frames}\label{sec-dual}

It is not difficult to see that 
the proof of the equivalences in 
Theorem \ref{prop1} hinges on the fact that $J$ happens to be a unitary operator 
with respect to $[\cdot,\cdot]$ and $[\cdot,\cdot]_J$. 
Theorem \ref{prop1} can be generalized by an application of 
the inverse mapping theorem as follows.
\begin{proposition}\label{pb}
Let $(\K_{1},[\cdot,\cdot]_{1})$ and $(\K_2,[\cdot,\cdot]_{2})$ 
be Krein spaces with chosen fundamental symmetries.  
Assume that $\{\eta_{x}^{1},\eta_{x}^{2},\ldots,\eta_{x}^{n}\}_{x\in\M}$  
is  a  continuous frame of rank $n\in\N$ 
for  $\K_{1}$ with respect to $(\M,\mathfrak{B},\nu)$ 
and with frame bounds $a\leq b$.  
Then, for any bijective linear operator $\boldsymbol{\A}\in\mathcal{B}(\K_{1},\K_{2})$, 
the family 
$\{\boldsymbol{\A}\eta_x^1,\boldsymbol{\A}\eta_x^2,\ldots,\boldsymbol{\A}\eta_x^n\}_{x\in\M}$ 
is a continuous frame of rank $n$ for $\K_{2}$ with frame bounds 
$a\hs \Vert\boldsymbol{\A}^{-1}\Vert^{-2} \leq\, b\hs \Vert\boldsymbol{\A}\Vert^{2}$.
\end{proposition}

The proof of the proposition is straightforward by observing that  
$|[\boldsymbol{\A}\eta_x^i,k]|=|[\eta_x^i,\boldsymbol{\A}^* k]|$, 
\,$\Vert\boldsymbol{\A}^*k\Vert_{J_1}\leq \Vert\boldsymbol{\A}\Vert\,\Vert k\Vert_{J_2}$ and 
\,$\Vert k\Vert_{J_2} \leq
\Vert\boldsymbol{\A}^{-1}\Vert\, \Vert\boldsymbol{\A}^*k\Vert_{J_1}$ 
for all $k\in\K_2$.  Proposition \ref{pb} motivates the following definition. 

\begin{definition}    \label{sf}
Let $(\K_{1},[\cdot,\cdot]_{1})$ and $(\K_2,[\cdot,\cdot]_{2})$ 
be Krein spaces. Two continuous frames  
$\{\eta_{x}^{1},\eta_{x}^{2},\ldots,\eta_{x}^{n}\}_{x\in\M}$ and 
$\{\theta_{x}^{1},\theta_{x}^{2},\ldots,\theta_{x}^{n}\}_{x\in\M}$ of rank $n\in\N$ 
with respect  to $(\M,\mathfrak{B},\nu)$ for the Krein spaces $\K_1$ and $\K_2$, respectively, 
are said to be \emph{similar} if there exists a bijective operator 
$\boldsymbol{\A}\in\mathcal{B}(\K_{1},\K_{2})$ 
such that $\boldsymbol{\A} \eta_{x}^{i}=\theta_{x}^ {i}$  for all $x\in\M$ and 
all $i=1,\ldots,n$. The continuous frames are called \emph{unitarily equivalent} if   
$\boldsymbol{\A}$ is a $J$-unitary operator. 
\end{definition}

Equations \eqref{desi:trans-Frame} and \eqref{Sinv} ensure that $\mathbf{S}$ 
belongs to $\mathcal{B}(\K)$ and has a bounded inverse. Therefore, the family 
$\{\mathbf{S}^{-1}\eta_{x}^{1},\mathbf{S}^{-1}\eta_{x}^{2},\ldots,
\mathbf{S}^{-1}\eta_{x}^{n}\}_{x\in\M}$ in the Frame Decomposition Theorem~\ref{Thm-trans-frame} 
yields a continuous frame of rank $n$ that is similar to 
$\{\eta_{x}^{1},\eta_{x}^{2},\ldots,\eta_{x}^{n}\}_{x\in\M}$. 

The next proposition provides a criterion for a pair of frames to be similar. 
It is an analogue of Deguang and Larson's result for Hilbert spaces \cite{DeguangLarson}.

\begin{proposition}\label{similarframesth} 
Let $(\K_{1},[\cdot,\cdot]_{1})$ and $(\K_2,[\cdot,\cdot]_{2})$ 
be Krein spaces. Two continuous frames  
$\{\eta_{x}^{1},\eta_{x}^{2},\ldots,\eta_{x}^{n}\}_{x\in\M}$ and 
$\{\theta_{x}^{1},\theta_{x}^{2},\ldots,\theta_{x}^{n}\}_{x\in\M}$ of rank $n\in\N$ 
with respect  to $(\M,\mathfrak{B},\nu)$ for the Krein spaces $\K_1$ and $\K_2$, respectively, 
are similar if and only if the corresponding analysis operators 
$\mathbf{T}_1:\K_1\ra \bigoplus_{i=1}^n\!\Elk(\M,\nu)$ and 
$\mathbf{T}_2:\K_2\ra \bigoplus_{i=1}^n\!\Elk(\M,\nu)$ 
have the same range.
\end{proposition}
\begin{proof}
The proof of the ``if'' part is the same as in \cite[Theorem 3.3]{Escobar-Esmeral-Ferrer} 
for $\K_{G}=\K_1$ and $\h_{\W}=\K_2$. It follows from %the fact that 
$[\boldsymbol{\A}\eta_{x}^{i},h]_2 = [\eta_{x}^{i}, \boldsymbol{\A}^*h]_1$ 
for all $h\in\K_2$, $i=1,\ldots,n$, and the bijectivity of 
$\boldsymbol{\A}\in\mathcal{B}(\K_{1},\K_{2})$. 

To prove the opposite direction, assume that 
$\operatorname{Rang}\mathbf{T}_{1}=\operatorname{Rang}\mathbf{T}_{2}=:V
\subset  \bigoplus_{i=1}^n\!\Elk(\M,\nu)$. 
Similarly to \eqref{Tfb}, we obtain from \eqref{ek}, \eqref{ipL2} and 
\eqref{oper-T[*]} that 
\[ \label{Tifb} 
a_i\hs \Vert k_i\Vert_{J_i}^2\ \leq\ \Vert \mathbf{T}_i 
k_i\Vert_{J_{\hsp\Elk}^\oplus}^2
\ \leq\  b_i\hs \Vert k_i\Vert_{J_i}^2  \quad 
\text{for \,all }\,k_i\hs\in\hs \K_i, 
\]
where $J_i$ denotes a chosen fundamental symmetry on $\K_i$ and 
$0<a_i\leq b_i$ are frame bounds for the frame in $\K_i$, \,$i=1,2$. 
From \eqref{Tifb}, we conclude that $\mathbf{T}_1$ and $\mathbf{T}_2$ 
are bounded, injective, and have closed range $V$ in the Hilbert space 
$\Big(\bigoplus_{i=1}^n\!\Elk(\M,\nu),[\cdot,\cdot]_{J_{\hsp\Elk}^\oplus}\Big)$. 
It follows that the operators $\mathbf{T}_1\in\mathcal{B}(\K_{1},V)$ and 
$\mathbf{T}_2\in\mathcal{B}(\K_{2},V)$ possess bounded inverses and the same 
holds for $\mathbf{T}_1^*\in\mathcal{B}(V,\K_{1})$ and 
$\mathbf{T}_2^*\in\mathcal{B}(V,\K_{2})$. 
Given $k_1\in\K_1$, set $k_2:=\mathbf{T}_2^{-1} \mathbf{T}_1 k_1$. 
Then obviously $\mathbf{T}_2k_2=\mathbf{T}_1 k_1$, which is equivalent to 
$[\theta_{x}^{i},k_2]_2 = [\eta_{x}^{i},k_1]_1\in \Elk(\M,|\nu|)$ 
for all $i=1,\ldots,n$. Recall from Section \ref{rk} that 
$V=\operatorname{Rang}\mathbf{T}_{1}=\operatorname{Rang}\mathbf{T}_{2}
\subset \bigoplus_{i=1}^n\!\Elk(\M,|\nu|)$ is a reproducing kernel Hilbert space, 
so $[\theta_{x}^{i},k_2]_2 = [\eta_{x}^{i},k_1]_1$ 
for all $x\in\M$ by \eqref{phiix}. 
With $k_2=\mathbf{T}_2^{-1} \mathbf{T}_1 k_1$, we obtain 
$$
[\eta_{x}^{i},k_1]_1 = [\theta_{x}^{i},\mathbf{T}_2^{-1} \mathbf{T}_1 k_1]_2 
= [\mathbf{T}_1^*\mathbf{T}_2^{*-1}\theta_{x}^{i},  k_1]_1 
\quad \text{for all }\,k_1 \in \K_1. 
$$
Since $(\K_{1},[\cdot,\cdot]_{1})$ is non-degenerate, 
it follows that $\boldsymbol{\A}\theta_{x}^{i}=\eta_{x}^{i}$ 
for all $x\in\M$ and $i=1,\ldots,n$, where 
$\boldsymbol{\A}:=\mathbf{T}_1^*\mathbf{T}_2^{*-1} 
\in\mathcal{B}(\K_{2},\K_{1})$ is a bounded bijective operator. 
This shows that the two continuous frames are similar. 
\end{proof}

\subsection{Dual and Parseval continuous frames}\label{sec-dual-similar}

This section is devoted to the description of \emph{dual} and \emph{Parseval} 
continuous frames and their relations to the the fundamental symmetry $J$, 
the frame operator $\mathbf{S}$ and coherent states. 
A similar analysis for discrete $J$-frames can be found in \cite{GMM}. 
We start with the definition of dual continuous frames.

\begin{definition}\label{def:dual-frame-|nu|}
Let $\{\eta_{x}^{1},\eta_{x}^{2},\ldots,\eta_{x}^{n}\}_{x\in\M}$ be a continuous frame of rank 
$n\in\N$ for the Krein space $(\K,[\cdot\hs,\cdot])$ with respect to $(\M,\mathfrak{B},\nu)$. 
A continuous frame $\{\hs \theta_x^1,\theta_x^2,\ldots,\theta_x^n\hs\}_{x\in\M}$ 
of rank $n$ 
\,for $(\K,[\cdot\hs,\cdot])$
with respect to $(\M,\mathfrak{B},\nu)$ is called a 
\emph{dual frame} of $\{\eta_{x}^{1},\eta_{x}^{2},\ldots,\eta_{x}^{n}\}_{x\in\M}$ 
(or sometimes a $\nu$-\emph{dual frame}) if 
$$
\uno_{\K} = \sum_{i=1}^{n}\int_{\M} |\theta_{x}^{i}]\hs[ \eta_{x}^{i}|\,\dd\nu(x).
$$	 
\end{definition}
Note that we used the measure $\nu$ instead of $|\nu|$ to emphasize the Krein space setting. 
Obviously, this can be compensated by multiplying $\theta_x^i$ by the value of the 
Radon-Nikodym derivative $j_{\Elk}$ at the point $x\in \M$. 
The following proposition is analogous to \cite[Proposition~3.8]{KEFER} 
in the discrete case. However, the Radon-Nikodym derivative $j_{\Elk}$ will appear 
in the formulas since the definition of the frame operator forced us to include 
the fundamental symmetry $J_{\hsp\Elk}^\oplus$, see Section \ref{impor-fund-symm}.  

\begin{proposition} \label{DF}
Let $(\K,[\cdot,\cdot])$ be a Krein space with fundamental symmetry $J$ and 
let $\{\eta_{x}^{1},\eta_{x}^{2},\ldots,\eta_{x}^{n}\}_{x\in\M}$ be a continuous frame 
of rank $n\in\N$ for $(\K,[\cdot,\cdot])$
with respect to $(\M,\mathfrak{B},\nu)$. With $\mathbf{S}$ denoting the corresponding 
frame operator given in \eqref{transformada-frame} and $j_{\Elk}$ denoting 
the Radon-Nikodym derivative of $\nu$ with respect to $|\nu|$, 
there are the following dual frames: 
\begin{enumerate}[i)]
\item  The continuous frame 
$\{j_{\Elk}\hsp(x)\hs\mathbf{S}^{-1}\eta_x^1,\,j_{\Elk}\hsp(x)\hs\mathbf{S}^{-1}\eta_x^2,\ldots,
j_{\Elk}\hsp(x)\hs\mathbf{S}^{-1}\eta_x^n\}_{x\in\M}$ is a dual frame of  
\,$\{\eta_{x}^{1},\eta_{x}^{2},\ldots,\eta_{x}^{n}\}_{x\in\M}$ 
\,for the Krein space $(\K,[\cdot,\cdot])$. 
\item \label{dii} The continuous frame 
$\{j_{\Elk}\hsp(x)\hs J\mathbf{S}^{-1}\eta_x^1,
\,j_{\Elk}\hsp(x)\hs J\mathbf{S}^{-1}\eta_x^2,\ldots,
j_{\Elk}\hsp(x)\hs J\mathbf{S}^{-1}\eta_x^n\}_{x\in\M}$ is a dual frame of 
\,$\{J\eta_x^1,J\eta_x^2,\ldots,J\eta_x^n\}_{x\in\M}$ 
\,for the Krein space $(\K,[\cdot,\cdot])$. 
\item \label{diii} The continuous frame 
$\{j_{\Elk}\hsp(x)\hs J\mathbf{S}^{-1}\eta_x^1,
\,j_{\Elk}\hsp(x)\hs J\mathbf{S}^{-1}\eta_x^2,\ldots,
j_{\Elk}\hsp(x)\hs J\mathbf{S}^{-1}\eta_x^n\}_{x\in\M}$ is a dual frame of 
\,$\{\eta_{x}^{1},\eta_{x}^{2},\ldots,\eta_{x}^{n}\}_{x\in\M}$ 
\,for the Hilbert space $\left(\K, [\cdot,\cdot]_{J}\right)$. 
\item \label{div} The continuous frame 
$\{j_{\Elk}\hsp(x)\hs\mathbf{S}^{-1}\eta_x^1,\,j_{\Elk}\hsp(x)\hs\mathbf{S}^{-1}\eta_x^2,\ldots,
j_{\Elk}\hsp(x)\hs\mathbf{S}^{-1}\eta_x^n\}_{x\in\M}$ is a dual frame of  
\,$\{J\eta_x^1,J\eta_x^2,\ldots,J\eta_x^n\}_{x\in\M}$ 
\,for the Hilbert space $\left(\K, [\cdot,\cdot]_{J}\right)$. 
\end{enumerate}
Furthermore, if $a\leq b$ are frame constants for 
$\{\eta_{x}^{1},\eta_{x}^{2},\ldots,\eta_{x}^{n}\}_{x\in\M}$, 
then $b^{-1}\leq a^{-1}$ are frame constants for all these dual frames. 
\end{proposition}
\begin{proof} 
Let $a\leq b$ be frame bounds for the frame 
$\{\eta_{x}^{1},\eta_{x}^{2},\ldots,\eta_{x}^{n}\}_{x\in\M}$. Since 
\begin{align*}
\sum_{i=1}^n\int_{\M}|[j_{\Elk}\hsp(x)\hs\mathbf{S}^{-1}\eta_x^i,k]|^2\hs \dd|\nu|(x)\, &=
\sum_{i=1}^n\int_{\M} [\mathbf{S}^{-1}k,\eta_x^i]\,[\eta_x^i,\mathbf{S}^{-1}k]\hs\dd|\nu|(x)\\
& =  [\mathbf{S}^{-1}k,\mathbf{S}\hs \mathbf{S}^{-1}k] = [k,J\mathbf{S}^{-1}k]_J
\end{align*} 
by \eqref{transformada-frame} and 
\,$b^{-1}\|k\|_J^2\leq [k,J\mathbf{S}^{-1}k]_J \leq a^{-1}\|k\|_J^2$ 
\,by the second relation in \eqref{Sinv}, it follows that 
$\{j_{\Elk}\hsp(x)\hs\mathbf{S}^{-1}\eta_x^1,\,j_{\Elk}\hsp(x)\hs\mathbf{S}^{-1}\eta_x^2,\ldots,
j_{\Elk}\hsp(x)\hs\mathbf{S}^{-1}\eta_x^n\}_{x\in\M}$ defines a continuous frame for 
$(\K,[\cdot,\cdot])$ with frame bounds $b^{-1}\leq a^{-1}$. Furthermore, 
from \eqref{transformada-frame} and Proposition~\ref{pf}.\ref{pfii}), we get 
$$
\uno_{\K} = \mathbf{S}^{-1} \mathbf{S} =\sum_{i=1}^{n}\int_{\M} 
|j_{\Elk}\hsp(x)\hs \mathbf{S}^{-1}\eta_{x}^{i}]\hs [\eta_{x}^{i}|\,\dd\nu(x), 
$$
so that it is indeed a dual frame of 
$\{\eta_{x}^{1},\eta_{x}^{2},\ldots,\eta_{x}^{n}\}_{x\in\M}$. 
Replacing $\mathbf{S}$ by $\mathbf{S}_0=J\mathbf{S}J$, $\mathbf{S}_1=\mathbf{S}J$ and 
$\mathbf{S}_2=J\mathbf{S}$ from Section \ref{cf} shows 
\,\ref{dii}), \ref{diii}) and \ref{div}), 
respectively. 
\end{proof}

The basic step in the proof of the last proposition was to multiply the frame operator 
by its inverse from the left and to apply Proposition \ref{pf}.\ref{pfii}). Likewise,  
if we multiply the frame operator from both sides 
with the inverse of a self-adjoint square root, 
then we will get a ``self-dual'' frame. However, this can only be done if  
a self-adjoint square root of the frame operator exists, for instance, by \eqref{Sinv}, 
for the positive Hilbert space operators $\mathbf{S}_1=\mathbf{S}J$ and 
$\mathbf{S}_2=J\mathbf{S}$. According to Theorem \ref{prop1}, the resulting 
frames for the associated Hilbert space will also yield frames for the Krein space. 
This observation is the starting point of the next proposition. 

\begin{proposition}\label{prop:raiz-frame-krein}
Let $(\K,[\cdot,\cdot])$ be a Krein space with fundamental symmetry $J$ 
and let $\{\eta_{x}^{1},\eta_{x}^{2},\ldots,\eta_{x}^{n}\}_{x\in\M}$ be a continuous frame 
of rank $n\in\N$ for $\K$  with respect to $(\M,\mathfrak{B},\nu)$. 
Then $\{(\mathbf{S}J)^{-1/2}\eta_x^1,\ldots,(\mathbf{S}J)^{-1/2}\eta_x^n\}_{x\in\M}$ and 
$\{(J\mathbf{S})^{-1/2}J\eta_x^1,\ldots,(J\mathbf{S})^{-1/2}J\eta_x^n\}_{x\in\M}$ 
define Parseval continuous frames of rank $n\in\N$ for $\K$, where 
$\mathbf{S}$ denotes the frame operator given in \eqref{transformada-frame} 
and $\K$ is viewed either as the Krein space $(\K,[\cdot,\cdot])$ 
or as the Hilbert space $(\K,[\cdot,\cdot]_J)$. 
\end{proposition}
\begin{proof}
Let $0<a\leq b$ be frame bounds for $\{\eta_{x}^{1},\ldots,\eta_{x}^{n}\}_{x\in\M}$ 
and thus, by Theorem \ref{prop1}, also for 
$\{J\eta_{x}^{1},\ldots,J\eta_{x}^{n}\}_{x\in\M}$. 
As is well-known from C*-algebra theory,  \eqref{Sinv} implies 
$$
0<b^{-1/2}\uno \leq (J\mathbf{S})^{-1/2}\leq a^{-1/2}\uno\quad 
\text{and} \quad 0<b^{-1/2}\uno \leq (\mathbf{S}J)^{-1/2}\leq a^{-1/2}\uno. 
$$
As a consequence, $(J\mathbf{S})^{\pm 1/2}$ and $(\mathbf{S}J)^{\pm 1/2}$ are bounded 
operators. Therefore, by Proposition~\ref{pb}, 
$\{(\mathbf{S}J)^{-1/2}\eta_x^1,\ldots,(\mathbf{S}J)^{-1/2}\eta_x^n\}_{x\in\M}$ and 
$\{(J\mathbf{S})^{-1/2}J\eta_x^1,\ldots,(J\mathbf{S})^{-1/2}J\eta_x^n\}_{x\in\M}$ 
define continuous frames for $\K$. 
Next, we compute for $k\in\K$ that 
\begin{align*}
&\sum_{i=1}^{n}\int_{\M}\left|[(\mathbf{S}J)^{-1/2}\eta_{x}^{i},k]\right|^{2}\dd |\nu|(x)=  
\sum_{i=1}^{n}\int_{\M}[J(\mathbf{S}J)^{-1/2}Jk, \eta_{x}^{i}]\, 
[\eta_{x}^{i},J(\mathbf{S}J)^{-1/2}Jk]\hs\dd |\nu|(x) \\
&=[J(\mathbf{S}J)^{-1/2}Jk, \mathbf{S}J (\mathbf{S}J)^{-1/2}Jk] 
= [Jk, (\mathbf{S}J)^{-1/2}(\mathbf{S}J)(\mathbf{S}J)^{-1/2}Jk]_J 
= [Jk,Jk]_J = \|k\|_J^2,  
\end{align*}
where we used $\big((\mathbf{S}J)^{-1/2}\big)^* 
= J\hs \big((\mathbf{S}J)^{-1/2}\big)^{[*]}J=J(\mathbf{S}J)^{-1/2}J$ 
in the first equality and \eqref{transformada-frame} in the second. 
This shows that $\{(\mathbf{S}J)^{-1/2}\eta_x^1,\ldots,(\mathbf{S}J)^{-1/2}\eta_x^n\}_{x\in\M}$ 
yields a Parseval continuous frame for $\K$, where, by Theorem \ref{prop1}, we may view $\K$ 
either as the Krein space $(\K,[\cdot,\cdot])$ or as the Hilbert space $(\K,[\cdot,\cdot]_J)$.

Likewise, since $\mathbf{S}_2=J\mathbf{S}$ functions as the frame operator 
for $\{J\eta_{x}^{1},\ldots,J\eta_{x}^{n}\}_{x\in\M}$ with respect 
to the Hilbert space $(\K,[\cdot,\cdot]_J)$, we obtain from \eqref{S2} that 
\begin{align*}
&\sum_{i=1}^{n}\int_{\M}\left|[(J\mathbf{S})^{-1/2}J\eta_{x}^{i},k]\right|^{2}\dd |\nu|(x)=  
\sum_{i=1}^{n}\int_{\M}[(J\mathbf{S})^{-1/2}Jk, J\eta_{x}^{i}]_J\, 
[J\eta_{x}^{i},(J\mathbf{S})^{-1/2}Jk]_J\hs\dd |\nu|(x) \\
&=[(J\mathbf{S})^{-1/2}Jk, J\mathbf{S} (J\mathbf{S})^{-1/2}Jk]_J
= [Jk, (J\mathbf{S})^{-1/2}(J\mathbf{S})(J\mathbf{S})^{-1/2}Jk]_J 
= [Jk,Jk]_J = \|k\|_J^2,  
\end{align*}
which shows that 
$\{(J\mathbf{S})^{-1/2}J\eta_x^1,\ldots,(J\mathbf{S})^{-1/2}J\eta_x^n\}_{x\in\M}$ 
defines a Parseval continuous frames for $\K$.  
\end{proof}

There are different notations of coherent states in the physics literature. 
For instance, Perelomov \cite[Section 2.1]{P} 
relates coherent states with unitary group representations, 
and Ali, Antoine and Gazeau \cite[Definition 5.4.5]{Ali-Antoine-Gazeau-2014} use the 
reproducing kernel property to define generalized coherent states. 
What both definitions have in common is that, in favorable situations 
(e.g.\ compact groups), a family of coherent states yields a resolution of the identity,   
commonly written as $\uno= \int_{X} |x\rangle \,\langle x|\,\dd \nu(x)$.  

\begin{example}[Coherent states in Krein spaces] \label{excs} 
Let $(\K,[\cdot,\cdot])$ be a Krein space with fundamental symmetry $J$ and 
let $\{\eta_{x}\}_{x\in\M}$ be a continuous frame  of rank $1$ for $(\K,[\cdot,\cdot])$ 
with respect to the \emph{positive} measure space $(\M,\mathfrak{B},\nu)$. 
Define $\sigma : \M \ra \K$ \,by $\sigma(x)\hsp :=\hsp (\mathbf{S}J)^{-1/2}\eta_{x}$. 
Then, by Proposition~\ref{prop:raiz-frame-krein},  $\{\sigma(x)\}_{x\in\M}$ 
yields a continuous Parseval frame for $\K$. 
Note that, for any Parseval frame with frame operator $S$, \eqref{Sinv} implies 
$S^{-1}J=\uno$, hence $S=J$. Thus \eqref{transformada-frame} gives 
$$
J =
\int_{\M}|\sigma(x)]\,[\sigma(x)|\,\dd\nu(x). 
$$ 
From this perspective, we may view the family $\{\sigma(x)\}_{x\in\M}$ 
as a system of coherent states for the Krein space $(\K,[\cdot,\cdot])$. 
\end{example}

In the light of Definition \ref{sf}, Proposition \ref{DF} states that any continuous 
frame is similar to a dual frame, 
and Proposition \ref{prop:raiz-frame-krein} shows that any continuous 
frame is similar to a Parseval frame. Moreover, combining Proposition \ref{r1} 
and Example \ref{excs}, we can say that any continuous frame 
with respect to a positive measure space determines a set 
of coherent states. 

The conclusion in Example \ref{excs} that the frame operator $\mathbf{S}$ 
of any Parseval frame satisfies $\mathbf{S}=J$ remains valid for 
continuous frames of rank $n$. Then \eqref{transformada-frame}  and 
Proposition \ref{pf}.\ref{pfii}) yield 
$$
\uno_{\K}  =\sum_{i=1}^{n}\int_{\M} 
|j_{\Elk}\hsp(x)\hs J\eta_{x}^{i}]\hs [\eta_{x}^{i}|\,\dd\nu(x). 
$$
From this, we conclude that, for any Parseval continuous frame 
$\{\eta_{x}^{1},\eta_{x}^{2},\ldots,\eta_{x}^{n}\}_{x\in\M}$, 
a dual frame is given by $\{j_{\Elk}\hsp(x)\hs J\eta_{x}^{1},
j_{\Elk}\hsp(x)\hs J\eta_{x}^{2},\ldots,j_{\Elk}\hsp(x)\hs J\eta_{x}^{n}\}_{x\in\M}$.

%----------------------
%---------
%\clearpage
\section{Continuous frames  in Krein spaces arising from a non-regular W-metric}
\label{Sec-W}

Each Krein space $(\K,[\cdot,\cdot])$ with fundamental symmetry $J$ can be given 
by considering the Hilbert $(\K,[\cdot,\cdot]_J)$ and setting 
$[h,k]:=[k,Jh]_J$ for all $k,h\in \K$. More generally, if 
$\mathbf{W}$ is a bounded self-adjoint operator on a Hilbert space $(\h,\ip{\cdot})$ 
with polar decomposition $\bW=J\hs |\bW|$ and 
such that $0\notin\spec(\mathbf{W})$, then 
\[ \label{ipW}
[f,g]:= \ip[f]{\mathbf{W}\hs g},\quad f,\, g \in\h, 
\] 
defines a non-degenerate inner product on $\h$.  
Since, by assumption,  $0\notin \spec(\mathbf{|W|})$, there exists $\epsilon>0$ 
such that $\epsilon\uno\leq |\bW| \leq \|\bW\|\hs \uno$. 
It follows that the norms $\|\cdot\|_J=\sqrt{\ip{|\bW|(\cdot)}}$ and 
$\|\cdot\|=\sqrt{\ip{\cdot}}$ are equivalent. Therefore, 
$(\h, [\cdot,\cdot])$ yields a Krein space with fundamental symmetry $J$ 
and fundamental decomposition $\h=\boldsymbol{\Pn}_+\h \oplus \boldsymbol{\Pn}_-\h$, 
where $\boldsymbol{\Pn}_\pm = \frac{1}{2}(\uno \pm J)$, 
see for instance \cite{Azizov}. 
Clearly, Krein or Hilbert spaces 
with the topology given by equivalent norms 
admit the same collection of (continuous) frames. 
Hence, by Proposition \ref{pb} and Theorem~\ref{prop1}, 
any continuous frame for $(\h,\ip{\cdot}) $ of rank $n\in\N$ yields one for 
$(\h, [\cdot\hs,\cdot])$ and vice versa.

As in \cite{KEFER}, the aim of our last section is to show 
how to transfer continuous frames for a Hilbert space $\h$ to 
a Krein space with a W-metric like the one defined in \eqref{ipW}, 
but with a possibly unbounded operator $\bW$ and allowing 
$0\in \spec(\bW)$. To begin, let 
$\mathbf{W}$ denote a self-adjoint operator with domain $\dom(\mathbf{W})\subset \h$,  
polar decomposition $\mathbf{W}= J\hs |\mathbf{W}|$, and integral representation 
$\mathbf{W}= \int \lambda \, \dd \mathbf{E}(\lambda)$, where   
$\mathbf{E}$ stands for the corresponding projection-valued measure on the 
Borel $\Sigma$-algebra $\mathfrak{B}(\R)$. 
We assume that $\ker(\mathbf{W})=\mathbf{E}(\{0\})=\{0\}$, which implies that  
$J$ is  a unitary self-adjoint operator.  
Analogous to  \eqref{ipW}, we define 
\[ \label{fg}
[f,g]:= \ip[f]{\mathbf{W}\hs g},\quad f,\, g \in \dom(\mathbf{W}). 
\] 
Then $\dom(\mathbf{W})$ becomes a decomposable non-degenerate inner product space 
with fundamental decomposition 
$\dom(\mathbf{W}) = \D_+\oplus \D_-$ and fundamental symmetry $J$, where 
\[
\D_+\hsp:=\hsp \mathbf{E}(0,\infty)\hs\dom(\bW),\ \ 
\D_-\hsp:=\hsp \mathbf{E}(-\infty,0)\hs\dom(\bW),\ \ 
J\hsp =\hsp \mathbf{E}(0,\infty)  \hsp-\hsp \mathbf{E}(-\infty,0). 
\]
Here, $\ker(\mathbf{W})=\{0\}$ is necessary since otherwise $\dom(\mathbf{W})$ would be degenerate. 
From $J^2 = \uno$, the polar decomposition $\mathbf{W}= J\hs |\mathbf{W}|$ 
and \eqref{fg}, it follows that 
\[  \label{Jip}
[f,g]_J = \ip[f]{|\mathbf{W}|\hs g},\quad f,\, g \in \dom(\mathbf{W}). 
\]
Taking the closure under the norm defined by $[\cdot\hs,\cdot]_J$ and 
extending  $J$ to the closure (without changing the notation), we obtain a Krein space 
$(\h_\bW, [\cdot\hs,\cdot])$ 
with fundamental symmetry $J$ and fundamental decomposition 
$\h_\bW = \h_+ \oplus \h_-$
such that $\D_+$ and $\D_-$ are dense in $\h_+$ and $\h_-$, respectively. 
The linear operator $\mathbf{W}$ in \eqref{ipW} and \eqref{fg} is called Gram operator. 
In case  $\bW\in \mathcal{B}(\h)$ and $0\notin\spec(\mathbf{W})$, 
the Gram operator $\bW$ and the Krein space $\left(\h_{\bW},[\cdot,\cdot]\right)$ 
are said to be regular.  

The next proposition summarizes some properties of $\h_{\bW}$. 
Its proof can be found in~\cite{KEFER}. 

\begin{proposition}\label{V3}
Let $\mathbf{W}\colon \dom(\mathbf{W}) \lra \h$ be a self-adjoint operator on 
the Hilbert space $\h$ such that $\ker(\mathbf{W})=\{0\}$. Then 
\begin{enumerate}[i)]
\item \label{Wi}
$\dom(\sqrt{|\mathbf{W}|}\hs)$ is complete in the norm $\|\cdot\|_J$ if and only if 
\,$0\notin \spec(\mathbf{W})$. 
In this case, 
%\item 
$\h_\bW$ can be identified with 
$\dom(\sqrt{|\mathbf{W}|})\subset \h$. 
%provided that $0\notin \spec(\mathbf{W})$. 
\medskip

\item \label{Wii}
$\sqrt{|\mathbf{W}|}\hs :\hs  \dom(\sqrt{|\mathbf{W}|}\hs)\hs \lra \hs  \h$ \,defines 
an isometric operator that 
admits an extension to a $J$-unitary operator 
\,$\mathbf{U}:=\overline{\phantom{X^+X^*}}\hspace{-24pt} {\sqrt{|\mathbf{W}|}}
\, : \,\h_\bW\,\lra\, \h$. 
\end{enumerate}
\end{proposition}

Proposition \ref{V3}.\ref{Wi}) shows that if $0\notin \spec(\mathbf{W})$, 
then we may set $\h_\bW:=\dom(\sqrt{|\mathbf{W}|}\hs)$,  
and if $0\in \spec(\mathbf{W})$, 
then the completion process will always require to 
add abstract elements to $\dom(\sqrt{|\mathbf{W}|})$.

\begin{example}\label{exam-2}
Let $(\M,\mathfrak{B},\mu)$ be measure space with a \emph{positive} and 
$\Sigma$-finite measure $\mu$. 
Given a measurable real function $\varphi: \M\ra \Real$, set 
$$
(\mathbf{W}_{\hsp \varphi} f)(p):= \varphi(p)\hs  f(p), \quad  
f\in \dom(\mathbf{W}_{\hsp \varphi}):=
\Big\{g\in 	\Elk(\M,\mu): \!\int_\M\! |g(p)|^2 |\varphi(p)|^2 \dd \mu(p)<\infty\Big\}, 
$$
and 
\[ \label{L2}
	[f,g]:= \ip[f]{\mathbf{W}_\varphi g}_{\Elk(\M,\mu)} 
	=\int_{\M}  \overline{f(p)} \hs g(p)\hs \varphi(p)\hs \dd\mu(p),\quad 
	f,g\in \dom(\mathbf{W}_{\hsp \varphi}). 
\]
If $0< \mathrm{ess\;inf}\hs  |\varphi | \leq \mathrm{ess\;sup}\hs |\varphi | <\infty$, 
then \eqref{L2} defines an (indefinite) inner product on 
$\h_{\bW_{\hsp \varphi}}=\dom(\mathbf{W}_{\hsp \varphi})=\Elk(\M,\mu)$
such that $(\Elk(\M,\mu), [\cdot\hs,\cdot])$ 
becomes a Krein space with fundamental symmetry $\je$ given by  
multiplication by the sign function $\mathrm{sign}(\varphi)$ of $\varphi$, 
i.e., $(\je f)(p):=\mathrm{sign}(\varphi)(p)\hs f(p)$. 
Note that this Krein space coincides with that in Example~\ref{exam-1} for 
$\nu:= \varphi\hs \mu$ and $|\nu| = |\varphi|\hs \mu$. 

Now assume that $ \mathrm{ess\;inf}\hs  |\varphi |>0$ and 
$\mathrm{ess\;sup}\hs |\varphi | = \infty$. In this case, 
$$
\h_{\bW_{\hsp \varphi}} 
=\Big\{g\in 	\Elk(\M,\mu): \!\int_\M\! |g(p)|^2 \hs |\varphi(p)|\hs \dd \mu(p)<\infty\Big\} 
= \dom(\mathbf{W}_{|\varphi|^{1/2}}) \subset \Elk(\M,\mu),  
$$
which enables us to view $\h_{\bW_{\hsp \varphi}}$ as a linear subspace of $\Elk(\M,\mu)$. 

If $ \mathrm{ess\;inf}\hs  |\varphi |=0$, we need to require 
$\mu(\{\varphi=0\})=0$ so that $\ker \mathbf{W}_{\hsp \varphi}=\{0\}$. Then 
$$
\h_{\bW_{\hsp \varphi}} 
=\Big\{\,g:\M\ra \Complex : \text{ measurable},\  
\int_\M\! |g(p)|^2 \hs |\varphi(p)|\hs \dd \mu(p)<\infty\,\Big\}, 
$$
but we cannot claim that $\h_{\bW_{\hsp \varphi}}\subset \Elk(\M,\mu)$ 
since $0\in \spec(\bW_{\hsp \varphi})$.   
On the contrary, whenever $\mathrm{ess\;sup}\hs |\varphi | <\infty$, 
we have $\Elk(\M,\mu) \subset \h_{\bW_{\hsp \varphi}}$.  

Nevertheless, in any case, the positive inner product $[\cdot,\cdot]_{\je}$ 
is given by 
$$
[ f,g]_{\je}= \int_{\M}  \overline{f(p)} \hs g(p)\hs |\varphi(p) |\hs \dd\mu(p),
\quad f,g\in \h_{\bW_{\hsp \varphi}}=\Elk(\M,|\varphi|\hs\mu).  
$$
\end{example}

In \cite[Proposition 4.1]{KEFER}, it has been shown that any discrete frame 
for a Hilbert space $\h$ cannot be a frame for a non-regular Krein space 
$(\h_\bW, [\cdot\hs,\cdot])$. The proof remains valid for continuous frames 
$\{\eta_{x}^{1},\eta_{x}^{2},\ldots,\eta_{x}^{n}\}_{x\in\M}$ of rank $n$ for 
$(\h,\ip{\cdot})$ and relies on the fact that 
$$
a\hs \|\bW k\|^2 \leq \sum_{i=1}^n\int_{\M}|\ip[\eta_x^i]{\bW k}|^2\hs\dd|\nu|(x) 
= \sum_{i=1}^n\int_{\M}|[\eta_x^i,k]|^2\hs\dd|\nu|(x) 
\leq b\hs \|\bW k\|^2 
\ \  \text{for all } k\in\h, 
$$
where $a\leq b$ are frame bounds for the continuous frame in $(\h,\ip{\cdot})$. 
Thus, if $0\in \spec(\bW)$ 
and $0\neq k\in \mathbf{E}(-\varepsilon,\varepsilon)\h$ 
for arbitrary small $\varepsilon>0$, 
%such that $\|k\|^2_J=  \ip[k]{|\mathbf{W}|\hs k}=1$,
%= \ip[\sqrt{|\mathbf{W}|}\hs k]{\sqrt{|\mathbf{W}|}\hs k} 
then $\bW k\neq 0$ and 
$$ 
 \|\bW k\|^2=  \ip[\!\sqrt{|\mathbf{W}|}\hs k]{|\bW|\sqrt{|\mathbf{W}|}\hs k} 
\leq \varepsilon \hs \ip[\!\sqrt{|\mathbf{W}|}\hs k]{\!\sqrt{|\mathbf{W}|}\hs k} 
=\varepsilon\hs \|k\|_J^2, 
$$
so that a lower frame bound does not exist. 
Similarly, if $\bW$ is an unbounded operator, we choose 
 $0\neq k\in \mathbf{E}\big(\hs\Real\hsp\setminus\hsp (-r,r)\hs\big)\h$ 
for arbitrary large  $r >0$ and obtain  
$$ 
 \|\bW k\|^2= \ip[\!\sqrt{|\mathbf{W}|}\hs k]{|\bW|\sqrt{|\mathbf{W}|}\hs k} 
\geq r \hs \ip[\!\sqrt{|\mathbf{W}|}\hs k]{\!\sqrt{|\mathbf{W}|}\hs k} 
=r\hs  \|k\|_J^2, 
$$
which proves that an upper frame bound cannot exist. 

These observations show that, for a non-regular $(\h_\bW, [\cdot\hs,\cdot])$, 
a continuous frame for $(\h,\ip{\cdot})$ will never be a continuous frame for 
$(\h_\bW, [\cdot\hs,\cdot])$. 
However, using the $J$-unitary operator $\mathbf{U}$ from Proposition \ref{V3}.\ref{Wii}) 
and applying Proposition \ref{pb} makes it possible to transfer 
continuous frames back and forth between $\h$ and $\h_{\bW}$. 
Note that a $J$-unitary operator does not change the frame bounds, 
and that $\mathbf{U}^{[*]} \!=\hsp\mathbf{U}^{-1}\! =\! \sqrt{|\bW|}^{-1}$ 
on \,$\dom(\!\sqrt{|\bW|}^{-1})$. 
We state the result in the following corollary.

\begin{corollary} \label{C}
Let \,$\mathbf{W}\colon \dom(\mathbf{W}) \ra \h$ \,be a self-adjoint operator on 
the Hilbert space $\h$ such that $\ker(\mathbf{W})=\{0\}$, 
and let \,$\mathbf{U}  : \h_\bW \ra  \h$ \,denote the $J$-unitary operator 
from Proposition \ref{V3}.\ref{Wii}). Then 
\begin{enumerate}[i)] 
\item \label{Ci}
$\{\eta_{x}^{1},\ldots,\eta_{x}^{n}\}_{x\in\M}$ is a continuous frame   
of rank $n\in\N$ for the Hilbert space $\h$ with respect to $(\M,\mathfrak{B},\nu)$ 
and with frame bounds $a\leq b$  if and only if 
$\{\mathbf{U}^{-1}\eta_{x}^{1},\ldots,\mathbf{U}^{-1}\eta_{x}^{n}\}_{x\in\M}$ 
is a continuous frame of rank $n$ for the Krein space $\h_\bW$ 
with respect to $(\M,\mathfrak{B},\nu)$ and with frame bounds $a\leq b$. 
\item \label{Cii}
If $\{\eta_{x}^{1},\ldots,\eta_{x}^{n}\}_{x\in\M}\subset \dom(\!\sqrt{|\bW|}^{-1})$ 
is a continuous frame of rank $n\in\N$ for the Hilbert space $\h$ 
with respect to $(\M,\mathfrak{B},\nu)$, then the continuous frame 
$\{\mathbf{U}^{-1}\eta_{x}^{1},\ldots,\mathbf{U}^{-1}\eta_{x}^{n}\}_{x\in\M}$ 
\,for the Krein space $\h_\bW$ is given by 
$\{\hsp\sqrt{|\mathbf{W}|}^{-1}\eta_{x}^{1},\ldots,
\!\sqrt{|\mathbf{W}|}^{-1}\eta_{x}^{n}\hs \}_{x\in\M}$.
\end{enumerate}
\end{corollary}

In practice, it might be difficult to find an explicit expression for 
$\sqrt{|\mathbf{W}|}^{-1}$. If, however, $\mathbf{W}$ is given by 
a multiplication operator $\mathbf{W}_\varphi$ on a function space as in 
Example \ref{exam-2}, then $\sqrt{|\mathbf{W}|}^{-1}$ can be expressed by 
$\mathbf{W}_{|\varphi|^{-1/2}}$. 
We illustrate this in our final example for $\h=\Elk(\R,\mu)$. 
In a certain sense, this is the generic case since 
any self-adjoint operator on a separable Hilbert space is 
unitarily equivalent to a direct sum of multiplication operators,  
see e.g.\ \cite[Theorem VII.3]{RS}. 
\begin{example}
Let $\mu$ be a positive and $\Sigma$-finite Borel measure on $\mathfrak{B}(\R)$ 
and set $\h:= \Elk(\R,\mu)$. 
For a measurable real function $\varphi: \R\ra \R$, consider as in Example \ref{exam-2} 
the operator 
$$
(\mathbf{W}_{\hsp \varphi} f)(t):= \varphi(t)\hs  f(t), \quad  
f\in \dom(\mathbf{W}_{\hsp \varphi}):=
\Big\{g\in 	\Elk(\R,\mu)\hs : \!\int_\R\! |g(t)|^2 |\varphi(t)|^2 
\dd \mu(t)<\infty\Big\}. 
$$
To ensure that $\ker \mathbf{W}_{\hsp \varphi}=\{0\}$, we 
require that $\mu(\{\varphi=0\})=0$. 
Then, as it has been shown in Example \ref{exam-2}, 
$$
\h_{\bW_{\hsp \varphi}} 
=\Big\{\,g:\R\ra \Complex : \text{ measurable},\  
\int_\R\! |g(t)|^2 \hs |\varphi(t)|\hs \dd \mu(t)<\infty\,\Big\}. 
$$
From this, it follows immediately that 
$\sqrt{|\varphi|}\, g\in \Elk(\R,\mu)$ for all $g\in\h_{\bW_{\hsp \varphi}}$. 
Therefore the $J$-unitary operator \,$\mathbf{U}$ in Proposition \ref{V3}.\ref{Wii}) 
and its inverse $\mathbf{U}^{-1}$ can be given by 
\begin{align*} 
&\qquad&
\mathbf{U}=\mathbf{W}_{|\varphi|^{1/2}} &\,:\,\h_{\bW_{\hsp \varphi}}\lra\Elk(\R,\mu), &
(\mathbf{W}_{|\varphi|^{1/2}} f)(t)&:= \sqrt{|\varphi(t)|}\hs f(t), &\qquad&\\ 
&\qquad&\mathbf{U}^{-1}=\mathbf{W}_{|\varphi|^{-1/2}}&\,:\,
\Elk(\R,\mu)\lra\h_{\bW_{\hsp \varphi}},&
(\mathbf{W}_{|\varphi|^{-1/2}} f)(t)&:=\mbox{$\frac{1}{\sqrt{|\varphi(t)|}}$}\hs f(t).
&\qquad& 
\end{align*}
Now, according to Corollary \ref{C}.\ref{Ci}),  
any continuous frame $\{\eta_{x}^{1},\eta_{x}^{2},\ldots,\eta_{x}^{n}\}_{x\in\M}$ 
of rank $n\in\N$ for $\Elk(\R,\mu)$ with respect to $(\M,\mathfrak{B},\nu)$ 
determines a continuous frame 
$$
\Big\{\efrac{1}{\sqrt{|\varphi|}} \hs\eta_{x}^{1},
\efrac{1}{\sqrt{|\varphi|}}\hs\eta_{x}^{2},\ldots,
\efrac{1}{\sqrt{|\varphi|}} \hs\eta_{x}^{n}\Big\}_{x\in\M}
\subset \h_{\bW_{\hsp \varphi}} 
$$ 
of rank $n$ for the Krein space $\h_{\bW_{\hsp \varphi}}$ 
with respect to $(\M,\mathfrak{B},\nu)$. 
Note that we did not assume that 
$\{\eta_{x}^{1},\eta_{x}^{2},\ldots,\eta_{x}^{n}\}_{x\in\M}
\subset \dom(\!\sqrt{|\mathbf{W}_{\hsp \varphi}|}^{-1})$
as it has been done in Corollary \ref{C}.\ref{Cii}).  
\end{example}

	%--------acknowledgment-----------
\subsection*{Acknowledgments} 
The second named author wishes to thank the Universidad de Caldas for financial support and hospitality.
 This is  part of the second  author’s project ``Elementos aproximadamente invertibles en C*-\'algebras y sus aplicaciones en teor\'ia de operadores''.
The third author acknowledges partial financial support from the 
%the Coordinaci\'on de la Investigaci\'on Cient\'ifica (UMSNH)
CIC-UMSNH %project "Grupos cu\'anticos y geometr\'ia no conmutativa"
and from the CONACyT project A1-S-46784.

	%%%-----------------
	%\clearpage

\end{document}